\documentclass[11pt]{amsart}

\usepackage{palatino} %1
\usepackage{hyperref}

\usepackage{amsmath,amssymb,amsthm,verbatim, amsfonts,amscd,flafter,epsf, epsfig,graphicx,verbatim,pinlabel,mathrsfs}
\usepackage[all]{xy}
\usepackage{epsf}
\usepackage[abs]{overpic}
\usepackage{epstopdf}

\usepackage{xcolor}
\usepackage{mathtools}
\usepackage{scalerel}

\definecolor{bred}{rgb}{1,0,.2}
\definecolor{blue}{rgb}{0,0,1}

\newtheorem{theorem}{Theorem}[section]
\newtheorem{lemma}[theorem]{Lemma}

\newtheorem{corollary}[theorem]{Corollary}
\newtheorem{proposition}[theorem]{Proposition}
\newtheorem{question}[theorem]{Question}
\newtheorem{theorem2}{Theorem}

\theoremstyle{definition}
\newtheorem{definition}[theorem]{Definition}

\newtheorem{remark}[theorem]{Remark}

\newtheorem{example}[theorem]{Example}

\def\tw{\textit{tw}}

\def\leave#1{{}}

\def\A{\mathcal{A}}

\def\X{\mathbf{X}}

\def\B{\mathcal{B}}
\def\H{\mathcal{H}}

\def\tb{{\mathit{tb}}}

   \title{Heegaard splittings and the tight Giroux Correspondence}
\author[Joan Licata]{Joan Licata}
\address{Mathematical Sciences Institute, Australian National University}
\email{joan.licata@anu.edu.au}

\author[Vera V\'ertesi]{Vera V\'ertesi}
\address{University of Vienna}
\email{vera.vertesi@univie.ac.at}

\begin{document}
\begin{abstract} 
This paper presents a new proof of the Giroux Correspondence for tight contact $3$-manifolds using techniques from Heegaard splittings and convex surface theory.  We introduce \emph{tight Heegaard splittings} of arbitrary contact $3$--manifolds; these generalise the Heegaard splittings naturally induced by an open book decomposition adapted to a contact structure on the underlying manifold.  Via a process called \emph{refinement}, any tight Heegaard splitting determines an open book, up to positive open book stabilisation.  This allows us to translate moves relating distinct tight Heegaard splittings into moves relating their associated open books.  We use this relationship to show that every Heegaard splitting of a contact $3$-manifold may be stabilised to  a Heegaard splitting induced by  a supporting open book decomposition.  Finally, we prove the tight Giroux Correspondence, showing that  any pair of open book decompositions supporting a fixed tight contact structure become isotopic after a sequence of positive open book stabilisations.  

\end{abstract}
\maketitle

\section{Introduction}

In this paper we prove the Giroux Correspondence for tight contact $3$--manifolds, showing that positive open book stabilisation suffices to relate any pair of open book decompositions adapted to a  tight contact structure.   Our methods lie firmly in the realm of contact topology, encoding equivalence classes of contact structures through topological, rather than geometric, data.   We hope the use of convex surface theory will make this proof accessible to a broad topological audience.

First studied as a purely topological object, an open book decomposition realises a $3$--manifold as a link with a fibered complement \cite{Alexander}.  Stallings  proposed surgery operations to relate distinct open books on a fixed $3$--manifold, and Harer proved the sufficiency of these moves \cite{Stallings, Harer}.  Later, Thurston-Winkelnkemper showed that an open book decomposition determines an equivalence class of contact structures on the underlying $3$--manifold, and Giroux famously extended this result, characterising the open books adapted to a fixed contact structure on a $3$--manifold as those related by a single one of Stallings's moves, \emph{positive stabilisation} \cite{TW, Giob}.  Recently, Breen-Honda-Huang have independently proven an analogous characterisation for contact manifolds in all odd dimensions \cite{BHH}.  

Our approach begins with another classical decomposition, a Heegaard splitting of a $3$--manifold.  Any open book decomposition determines a Heegaard splitting, but not all Heegaard splittings can be realised through this process.  For example, every Heegaard splitting induced by an open book has Hempel distance less than or equal to $2$, while there exist Heegaard splittings with arbitrarily large distance; further comparisons between these may be found in \cite{Rub}.   

Work of Torisu allows us to identify when a Heegaard splitting of $(M,\xi)$ is induced by an open book supporting $\xi$, and we call such splittings \emph{convex Heegaard splittings} \cite{Torisu}.  Given an arbitrary Heegaard splitting of a contact $3$--manifold $(M, \xi)$, we show how to stabilise it to a convex splitting, and we note that one may recover an open book decomposition supporting $\xi$ from any convex splitting of $M$.   This provides a new and accessible proof of the following result:

\begin{theorem2}\label{thm:ex}
Any contact $3$--manifold admits a compatible open book decomposition. 
\end{theorem2}

 We define a more general notion of a Heegaard splitting compatible with a contact structure:

\begin{definition} A Heegaard splitting $(\Sigma, U, V)$ of a contact $3$--manifold $(M,\xi)$ is \emph{tight} if $\Sigma$ is convex and $\xi$ restricts tightly to each handlebody $U, V$.
\end{definition}

A convex Heegaard splitting is necessarily tight, but not every tight Heegaard splitting is induced by an open book \cite{Torisu}.  Nevertheless, we may construct an open book from any tight splitting.  We introduce a process called \textit{refinement} which involves stabilising a tight Heegaard splitting to a convex Heegaard splitting.  As noted above, this convex splitting determines an open book decomposition of the original manifold.  Refinement involves many choices -- a Heegaard diagram for the splitting, sets of convex compressing discs inducing the diagram, properly embedded arcs on said discs -- and different choices will produce different convex splittings.

\begin{figure}[h]
\begin{center}
\includegraphics[scale=.8]{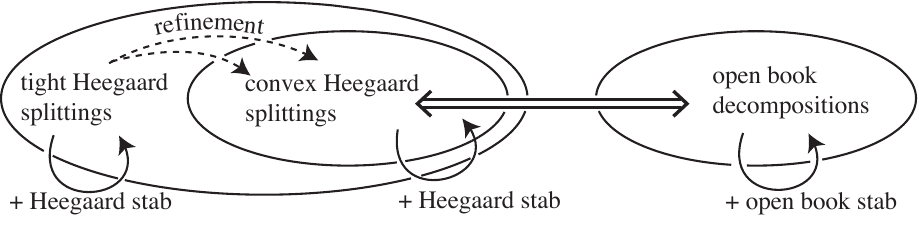}
\label{fig:intro2}
\end{center}
\end{figure}

Nevertheless, the associated open books are closely related: 

\begin{theorem}\label{thm:posstabclass} Suppose $(\Sigma, U, V)$ is a tight Heegaard splitting of $(M, \xi)$ and suppose that $(B, \pi)$ and $(B', \pi')$ are open book decompositions determined by refinements of $(\Sigma, U, V)$.  Then $(B, \pi)$ and $(B', \pi')$ admit a common positive open book stabilisation.  
\end{theorem}

When the original contact  $3$--manifold is tight, any convex Heegaard surface defines a tight splitting, so Theorem~\ref{thm:posstabclass} implies that the Heegaard splitting determines a  positive stabilisation class of open books; this answers a question of Rubinstein \cite{Rub}.  Torisu showed that every contact manifold -- whether tight or overtwisted -- admits a tight Heegaard splitting, and in each case we again recover a positive stabilisation class of open book decompositions \cite{Torisu}.

We also identify a move between tight Heegaard splittings (\textit{positive Heegaard stabilisation}) with the property that the refinements of splittings with a common positive Heeagard stabilisation yield open books with a common positive open book stabilisation.  This observation suggests that tight Heegaard splittings are worth studying in their own right, as opposed to  merely as a route to open books. We have included some  questions for further investigation in the next section.

 In this paper, however, the primary motivation for studying tight Heegaard splittings is a proof of the following theorem:

\begin{theorem2}\label{thm:GC}[Tight Giroux Correspondence] If two open book decompositions $(B,\pi)$ and $(B',\pi')$ of $M$ support isotopic tight contact structures, then they are related by a sequence of positive stabilisations and destabilisations. 
\end{theorem2}

Our proof begins by considering the pair of convex Heegaard splittings induced by a pair of open books for $M$ which support a fixed contact structure $\xi$.  The  Reidmeister-Singer Theorem asserts that any pair of Heegaard splittings for $M$ will become isotopic after sufficiently many Heegaard splitting stabilisations of each.  We may choose these stabilisations to correspond to stabilisations of the associated open books, leading us to study a pair of open book decompositions that induce isotopic Heegaard splittings.  The isotopy discretisation argument of Colin (stated as Theorem~\ref{thm:isdisc}) decomposes smooth isotopies of convex splitting surfaces as a sequence of convex isotopies and bypass attachments.  This produces a sequence of convex Heegaard surfaces related to each other by bypass attachments, and by the work above,  we may associate an open book decomposition supporting $\xi$ to each of these.  Finally, we show that these open books are related by positive open book stabilisation.

\subsection{Reading guide and open questions}
Section~\ref{sec:background1} provides background for the paper, stating well known technical results  in the generality that will be useful later on.   It also introduces  notation for the upcoming sections, but readers familiar with convex surface theory  are advised to restrict their attention to the discussion of bypasses in Section~\ref{sec:bypass}.  Section~\ref{sec:background2} discusses open books and Heegaard splittings in the context of contact geometry.  The first part may also be viewed as a background section, but one that establishes the perspective used in the main results that follow.   Section~\ref{sec:ob} introduces positive stabilisation for Heegaard splittings and presents some essential technical results. Section~\ref{sec:ex} presents a short proof that every contact $3$--manifold admits a supporting open book; this result is an application of ideas from the previous section, but is not used later on.  The technical heart of the paper lies in Section~\ref{sec:tighths}, where we define the refinement of a tight Heegaard splitting and show that the positive stabilisation class of the associated open book is well defined.   With these tools in hand, the proof of the tight Giroux Correspondence is short, and the final section discusses the potential for and obstacles to extending this approach to overtwisted manifolds.

Throughout the paper, we require all Heegaard surfaces to be convex, but we consider a variety of compatibility conditions between the contact structure and the handlebodies.  The strictest relationship is seen in the convex Heegaard splittings (Section~\ref{sec:obch}) directly induced by a supporting open book.  At the other extreme,  Section~\ref{sec:ex}  constructs an open book from an arbitrary Heegaard splitting, at a cost of greatly increasing the genus of the splitting via contact $1$-handle addition.  Section~\ref{sec:tighths} considers a middle ground, stabilising tight Heegaard splittings (Definition~\ref{def:tighths}) to yield convex splittings.  Rubinstein notes that in the Heegaard genus $2$ case, ``open book decompositions are nearly always more complicated than minimal Heegaard splittings'', and \"Ozba\u{g}ci shows the analogous statement in the genus $1$ case \cite{Rub}, \cite{Ozb} .  Although this is a topological observation rather than a contact one, it is consistent with the fact that we use stabilisation to transform an arbitrary Heegaard splitting to a convex splitting.  However, stabilisation destroys other information carried by the splitting; for example,   a Heegaard splitting of distance at least $3$ implies that the underlying manifold is hyperbolic.  

Below, we include some additional questions exploring the relationship between contact structures and Heegaard splittings.

\begin{question} Is there a bound on the distance of a tight Heegaard splitting?
\end{question}

\begin{question} Given a convex Heegaard surface in $(M, \xi)$, what is the minimal number of stabilisations required to make the splitting contact?
\end{question}

As indicated above, we prove the hard direction of the Giroux Correspondence only for tight contact $3$--manifolds.  Nevertheless, many of the notions developed in the paper apply equally well to tight splittings of overtwisted manifolds, leading to other natural questions:

\begin{question} For an overtwisted contact structure $\xi$ on $M$, what is the minimal genus of a tight Heegaard splitting?
\end{question} 

\begin{question} For an overtwisted contact manifold, can the Heegaard genus and the minimal genus of a tight splitting be arbitrarily far apart?
\end{question}

\subsection*{Acknowledgment} This material is based in part upon work supported by the National Science Foundation under Grant No. DMS-1929284 while the authors were in residence at the Institute for Computational and Experimental Research in Mathematics in Providence, RI, during the Braids program.  The first author received support from the Australian National University's Outside Studies Program and the second author was  supported by the FWF grant ``Cut and Paste Methods in Low Dimensional Topology'' P 34318.  The second author would like to thank the  the Erd\"{o}s Center  for providing a friendly and calm environment.  Both authors appreciate the helpful feedback provided by the referee.

\section{Background}\label{sec:background1} 
This section provides technical background that will be relied upon in the rest of the paper.   We assume familiarity with contact structures, Legendrian knots, and the Thurston-Bennequin number; readers unfamiliar with these are directed to \cite{Geiges}.  We have subsections on several topics:  \ref{sec:conv} convex surfaces, \ref{sec:handles} contact handles, and \ref{sec:bypass} bypasses.  The organisation is intended to be as transparent as possible for readers who wish to move to Section~\ref{sec:background2} or Section~\ref{sec:ex} and refer back as needed.  Henceforth, $M$ is always assumed to be an oriented $3$--manifold.

\subsection{Convex surfaces}\label{sec:conv} 
Here we give a brief introduction to the essentials of convex surface theory; for a more thorough treatment, the reader is referred to \cite{Ma}. 
Given   a contact structure $\xi=\ker{\alpha}$ on $M$, a surface $\Sigma$ embedded   in $M$ is \emph{convex} if there is a \emph{contact vector field} $X$ (i.e., a vector field whose flow preserves $\xi$) transverse to $\Sigma$.  If $\partial \Sigma\neq \emptyset$, we require $\partial \Sigma$ to be Legendrian.
Using the transverse direction given by $X$, one can build a neighbourhood  $\nu(\Sigma)\cong (\Sigma\times I,\xi\vert_{\nu(\Sigma)})$ with an $I$-invariant contact structure. (When $\Sigma \subset \partial M$, the transverse vector field yields an $I$-invariant half-neighbourhood.). The existence of such a (half-)neighbourhood is an alternative criterion for the convexity of $\Sigma$.

The locus of points where $\alpha(X)=0$ is a 1-dimensional submanifold $\Gamma_\Sigma$ called the \emph{dividing curve}.  The dividing curve separates $\Sigma$ into two submanifolds $\Sigma_\pm:=\{x: \pm \alpha(X)>0\}$.  Any embedded surface can be made convex via a $C^\infty$-isotopy, and an isotopy that keeps $\Sigma$ convex is called a \emph{convex isotopy}. For a given convex surface, the choice of transverse convex vector field is not unique, but different choices will yield dividing sets that differ only by  isotopy on $\Sigma$.

The dividing curve on a convex surface can be used to compute the relative twisting of the contact planes along Legendrian curves.  If $C$ is Legendrian on the convex surface $\Sigma$, then 
\[\textit{tw}_C(\xi,T\Sigma)=-\frac12|\Gamma_\Sigma\cap C|,\] 
where $\textit{tw}_C(\xi,T\Sigma)$ denotes the relative twisting of the plane fields $\xi\vert_C$ and $T\Sigma\vert_C$ along $C$.

 Let $L$ be a boundary component of $\Sigma$ and take a standard neighbourhood $(\nu(L),\xi\vert_{\nu(L)})$ contactomorphic to 
\[\left(S^1\times D^2,\xi=\ker(\sin(n\vartheta)\mathit{d}x+\cos(n\vartheta)\mathit{d}y)\right),\]
where $\vartheta$ is the coordinate parameterising $L$ and $x,y$ are coordinates on $D^2$. 
The boundary component $L$ is in \emph{standard position} if $\Sigma$ restricts to  $\{y=0,x\ge0\}$ in this model. 

As shown by Kanda \cite{Kanda}, boundary of a surface can be isotoped to standard position by a $C^0$-isotopy supported in $\nu(L)$ if and only if $\tw_L(\xi,T\Sigma)\le 0$. Once the boundary is in standard position, the interior of the surface can be isotoped to a convex position via a $C^\infty$-isotopy that fixes a neighbourhood of the boundary.

We can similarly measure relative twisting of Legendrian arcs with endpoints on the dividing set:
Let $\Sigma$ be a convex surface and $C$ a Legendrian arc with endpoints on $\Gamma_\Sigma$.  Then 
\[\textit{tw}_C(\xi,T\Sigma)=-\frac12|\Gamma_\Sigma\cap C|,\] 
where the endpoints of $C$ on $\Gamma_\Sigma$ are each counted with multiplicity $\frac{1}{2}$ on the right-hand side.

Consider two convex surfaces $\Sigma$ and $\Sigma'$ that intersect each other transversely along a closed Legendrian curve. 
The surfaces $\Sigma$ and $\Sigma'$ intersect \emph{standardly} along a component $L$ of $\Sigma\cap \Sigma'$ 
if, in the standard neighbourhood $\nu(L)$, we have \[\Sigma\cap\nu(L)=\{y=0\}\qquad \text{ and }\qquad \Sigma'\cap\nu(L)=\{x=0\}.\] If $L$ is a boundary component of $\Sigma$ and $\Sigma'$, then we require 
\[\Sigma\cap\nu(L)=\{y=0,x\ge0\}\text{ and }\Sigma'\cap\nu(L)=\{x=0,y\ge0\}.\] 

Again, Kanda \cite{Kanda} showed  that if $\tw_L(\xi,T\Sigma)=\tw_L(\xi,T\Sigma')\le 0$, then standard intersection can be achieved  -- in a slightly smaller neighbourhood of $L$ -- by a $C^0$-isotopy supported in $\nu(L)$ that keeps $L$ fixed and both $\Sigma$ and $\Sigma'$ convex. 
Once $\Sigma$ and $\Sigma'$ intersect  standardly, the intersection points $\Gamma_\Sigma\cap L$ and  $\Gamma_{\Sigma'} \cap L$ alternate along $L$.

The union of convex surfaces with Legendrian boundary gives a \emph{piecewise convex surface}. More precisely,  this is a surface $\Sigma=\cup\Sigma_i$ where each $\Sigma_i$ is convex with Legendrian boundary and such that the following hold:
\begin{enumerate}
\item for distinct $i,j,k$, $\Sigma_i\cap\Sigma_j\cap\Sigma_k=\emptyset$;
\item $\Sigma_i\cap\Sigma_j=\partial \Sigma_i\cap\partial \Sigma_j$ is Legendrian, and at each component of a double intersection,  the surfaces $\Sigma_i$ and $\Sigma_j$ intersect standardly.
\end{enumerate}

We can smooth a piecewise convex surface along any component of the double intersections. 

\begin{figure}[h]
\labellist
\small\hair 2pt
\pinlabel $\Sigma'$ [l] at 85 30
\pinlabel $\Sigma$ [l] at 6 10
\pinlabel {$L$} [tr] at 65 16
\pinlabel $\Sigma''$ [l] at 210 30
\endlabellist

\begin{center}
\includegraphics[scale=1.2]{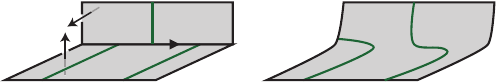}
\caption{ Left: Two convex surfaces meeting along a Legendrian curve.  Right: After smoothing, the new convex surface.} 
\label{fig:smoothing} 
\end{center}
\end{figure}
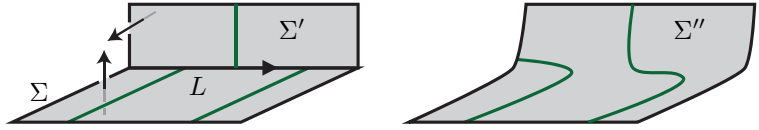

\begin{lemma}[Edge rounding]\label{lem:edgerounding} 
Let $\Sigma$ and $\Sigma'$ be convex surfaces with a standard intersection along some Legendrian curve $L$ which is a boundary component of $\Sigma$ and $\Sigma'$ and oriented as a boundary of $\Sigma'$. 
Then we may form a smooth convex surface $\Sigma''$ by replacing $\Sigma$ and $\Sigma'$ in $\nu(L)$ so that $\Gamma_{\Sigma''}$ restricts to $\Gamma_\Sigma$ and $\Gamma_{\Sigma'}$ away from $\nu(L)$ and connects each endpoint of  $\Gamma_{\Sigma'}$ to the next endpoint of  $\Gamma_{\Sigma}$ along the oriented curve $L$, as in Figure \ref{fig:smoothing}.

\end{lemma}

Note that in the local model $\nu(L)$, the surfaces $\{y=\varepsilon,x\ge\varepsilon\}$ and $\{x=\varepsilon,y\ge\varepsilon\}$ are convex and intersect standardly at $\{x=y=\varepsilon\}$.  This observation allows us to construct a tubular neighbourhood $\nu(\Sigma)$ for piecewise convex surfaces which is foliated by copies of the piecewise convex surface.  A \emph{convex isotopy} of a piecewise convex surface is an isotopy that keeps it piecewise convex at all times. A piecewise convex surface is \emph{closed} if each boundary component of the $\Sigma_i$ is in a double intersection. In particular, closed convex surfaces are  closed piecewise convex.

One can also introduce corners along Legendrian simple closed curves on a convex surface $\Sigma$, which will be useful for gluing contact manifolds with boundary.  Consider two $C^0$ convex isotopies of $\Sigma$ in a standard neighbourhood $\nu(L)$ of $L$ so that the images $\Sigma'$ and $\Sigma''$ in $\nu(L)$ meet transversely along $L$.  Form the new cornered surface  by cutting $\Sigma'$ and $\Sigma''$ along $L$ and then gluing a component of one to a component of the other inside $\nu(L)$.  There are two choices, depending on the preferred type of corner, and edge rounding either of the resulting cornered surfaces returns a smooth surface convexly isotopic to the original $\Sigma$.

The dividing curve determines the contact structure in an $I$-invariant neighbourhood of a convex surface.  Although often presented as a statement about what characteristic foliations can appear in an $I$-invariant neighbourhood, this important result can be rephrased as follows: 

\begin{theorem}\cite{Gi,Kanda}\label{thm:convex}
Let $\iota: \Sigma \hookrightarrow (M, \xi)$ and $\iota': \Sigma \hookrightarrow (M', \xi')$ be embeddings such that   $\iota^{-1}(\Gamma_{\iota(\Sigma)})={\iota'}^{-1}(\Gamma_{{\iota'}(\Sigma)})$. Then for any neighbourhood $N$ of $\iota(\Sigma)$ with $I$-invariant contact structure $\xi\vert_N$, there exists a neighbourhood $N'$ of $\iota'(\Sigma)$ and a contact embedding  $f:(N',\xi'\vert_{N'})\hookrightarrow (N,\xi\vert_{N})$. Moreover, the image  $f(\iota'(\Sigma))$ is convex isotopic to $\iota(\Sigma)$ via an isotopy that fixes $\Gamma_{\iota(\Sigma)}$ and is transverse to the $I$-direction.

If $\Sigma$ has non-empty boundary, then we assume that both $\iota(\Sigma)$ and $\iota'(\Sigma)$ are standard near $\iota(\partial \Sigma)$ and $\iota'(\partial \Sigma)$. In this case, the above isotopy fixes a neighbourhood of the boundary of $\iota(\Sigma)$. 
\end{theorem}

\begin{definition}
Two contact manifolds $(M,\xi)$ and $(M',\xi')$ with piecewise convex boundary are \emph{weakly contactomorphic} if there is a contact embedding $\iota\colon M\hookrightarrow M'$  such that $\iota(\partial M)$ is convex isotopic to $\partial M'$.
\end{definition}

Weak contactomorphism is an equivalence relation. It is clearly reflexive and transitive.  To establish symmetry, assume that $(M,\xi)$ embeds into $(M',\xi')$. Extend $(M',\xi')$ by an $I$-invariant neighbourhood of $\partial M'$ to obtain a weakly contactomorphic contact manifold $(M'',\xi'')$.  Then by Theorem~\ref{thm:convex}, one can find another copy $\Sigma$ of $\partial M$ in this extension; truncate $M''$ at $\Sigma$ to get a contact manifold weakly contact isotopic to $(M,\xi)$ that contains $(M',\xi')$.

\begin{definition}

Two embedded codimension-$0$ submanifolds $N_0,N_1\subset (M,\xi)$ with piecewise convex boundary are \emph{weakly (contact) isotopic} if there is an isotopy $N_s$ between them such that $\partial N_s$ is piecewise convex throughout. 
\end{definition}

The advantage of this notion is that it allows us to disregard the specific characteristic foliation on the boundary of $N_i$ and  concentrate only on $\Gamma_{\partial {{N_i}}}$.

\begin{theorem}[Giroux's Criterion]\label{thm:gircrit}
A convex surface $\Sigma$ in a contact manifold $(M,\xi)$ has a tight neighbourhood if and only if 
\begin{enumerate}
\item $\Sigma$ is a sphere and $\Gamma_\Sigma$ is connected; or
\item $\Sigma$ is not a sphere and $\Gamma_\Sigma$ has no closed contractible components.
\end{enumerate}
\end{theorem}

Suppose $(M,\xi)$ has a piecewise convex boundary, and let $L$ be a Legendrian corner between the convex pieces $\Sigma$ and $\Sigma'$. If $\Sigma=\{y=0,x\ge0\}$ and $\Sigma'= \{x=0,y\ge0\}$, $M=\{x,y\ge 0\}$ in the local model near $L$, then one can simply round $\Sigma$ and $\Sigma'$ at $L$ inside $M$, and declare the submanifold with boundary $\Sigma''$ to be the \emph{rounded} $M$.
On the other hand, suppose $M=\{y\le 0 \text{ or } x\le 0\}$ in the local model. Then we first take a parallel copy of $\Sigma$ inside $M$,  round so that $\Sigma''$ does not touch $\partial M$, and declare the manifold bounded by $\Sigma''$ to be the rounded $M$. This operation is independent of the choices made,  up to weak isotopy.

\begin{theorem}\cite{El}\label{thm:Eliashberg}
Let $B$ be a topological $3$-ball with piecewise convex boundary.  Up to weak contactomorphism on $B$ after rounding, there is a unique tight contact structure inducing a connected dividing set on the rounded ball.  \end{theorem}

There exists a large topological class of curves on a convex surface which have Legendrian representatives.  Quoting \cite{HKM}, an embedded graph $C$ on a convex surface $S$ is \textit{non-isolating} if $C$ is transverse to $\Gamma_S$; the univalent vertices of $C$ (and no others) lie on $\Gamma_S$; and every component of $S\setminus (\Gamma_S\cup C)$ has a boundary component which intersects $\Gamma_S$. 
\begin{theorem}[Legendrian Realisation Principle]\label{thm:lerp}

Given any non-isolating graph $C$ on a convex surface $S$, there exists a convex isotopy $\phi_s, s\in [0,1]$ such that the following hold:
\begin{enumerate}
\item $\phi_0=\text{id}$ and $\phi_s|_{\Gamma_S}=\text{id}$;
\item $\phi_1(\Gamma_S)=\Gamma_{\phi_1(S)}$; and
\item $\phi_1(C)$ is Legendrian.
\end{enumerate}
\end{theorem} 

\begin{proposition}[Partial gluing]\label{prop:glue}

Let  $(M^R,\xi^R)$ and $(M^L,\xi^L)$  be  contact 3-manifolds with piecewise convex boundary,  and suppose that $\varphi\colon \Sigma^R \to \Sigma^L$  is a diffeomorphism that identifies a pair of convex components of $\partial M^L$ and $\partial M^R$ and carries $\Gamma_{\Sigma^R}$ to $\Gamma_{\Sigma^L}$. Then up to weak contactomorphism, there is a unique contact structure $\xi:=\xi^L\cup\xi^R$ on $M=M^R\cup_{\varphi} M^L$ with a piecewise convex boundary 
 that restricts (again, up to weak contactomorphism) to $M^R$ and $M^L$ as $\xi^R$ and $\xi^L$. 
\end{proposition}

\begin{proof}

Constructing the glued-up contact structure follows the procedure of smooth gluing; see, for example, Section 2.7 in \cite{wall}. First, extend $M^R$ and $M^L$ along $\Sigma^R$ and $\Sigma^L$ by $I$-invariant contact structures  $\Sigma^R\times I$ and $\Sigma^L\times I$, respectively.  Use Theorem~\ref{thm:convex} to find a neighbourhood of $\Sigma^R$ in the extension of $M^L$. Truncate $M^L$ along the identified copy of $\Sigma_R$ and identify the neighbourhoods of $\Sigma^R$ in the two manifolds. 
The proof of uniqueness also follows the smooth argument; manifolds $M_1$ and $M_2$  resulting from different gluing choices are related by a contactomorphism with the following property.  In each $M_i$, there exists a neighbourhood $\nu_i$ of the image of $\Sigma_L$ which is weakly contact isotopic to an $I$-invariant neighbourhood.  There is a contactomorphism taking $M_1$ to $M_2$ that maps $\nu_1$ to $\nu_2$ and restricts to the identity away from this neighbourhood.
\end{proof}

\subsection{Contact handles}\label{sec:handles}
\begin{figure}[h]
\labellist
\small\hair 2pt
\pinlabel $L$ [l] at 127 18
\endlabellist

\begin{center}
\includegraphics[scale=1.2]{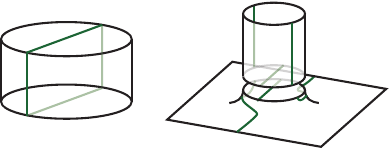}
\caption{ Left: a contact $1$- or $2$-handle. Right: Adding a corner along the closed Legendrian $L$ allows us to smoothly attach a contact $1$-handle.} 
\label{fig:handle}
\end{center}
\end{figure}

The basic building blocks for contact manifolds are contact handles.  These were first introduced by Giroux  \cite{Gi}, but in this paper we find it convenient to use a reformulation by \"Ozba\u{g}ci \cite{Oz} phrased in the language of convex surfaces.
Since every $3$--dimensional $k$--handle is topologically a ball, contact handles are simple to describe up to weak contactomorphism using Theorem \ref{thm:Eliashberg}: there is a unique tight contact structure on $D^3$ with smooth convex boundary and connected dividing set. This is the model for a \emph{contact 0-handle} $(h^0,\zeta^0)$ and a \emph{contact 3-handle} $(h^3,\zeta^3)$. Similarly, there is a unique tight contact structure on $D^1\times D^2$ with dividing curve $\Gamma$ as in Figure \ref{fig:handle};  this is the model for both a \emph{contact 1-handle} $(h^1,\zeta^1)$ and a \emph{contact 2-handle} $(h^2,\zeta^2)$. Coordinate models for contact handles of each index may be found in \cite{Oz}. 

As usual, contact $0$-handles are attached to the empty set, but the attaching data is important for higher-index contact handles.
Given a contact cobordism $(W,\xi)$ with convex boundary $\partial_+W\cup-\partial_-W$, let $\varphi^1\colon \partial D^1\times D^2\hookrightarrow \partial_+W$ be a diffeomorphism such that each $D^2$ component of the image of $\varphi^1$ is intersected in an arc by the dividing curve $\Gamma_{\partial_+W}$. Choose a representative of $\xi$ so that $L=\varphi^1(\partial D^1\times \partial D^2)$ is Legendrian. Introduce a corner along $L$ as described in Section~\ref{sec:conv};  $\varphi^1$ is still a map into the  newly piecewise convex surface, which we  still denote by $\partial_+W$. Using Proposition \ref{prop:glue}, glue $h^1$ onto $W$ to obtain a new contact manifold $(W\cup h^1,\xi\cup\zeta^1)$.  By construction, this manifold already has smooth boundary. This is a \emph{contact 1-handle attachment}.

To attach a \emph{contact 2-handle}, start with a diffeomorphism $\varphi^2\colon \partial D^2\times D^1\hookrightarrow \partial_+W$ with the property that the two arcs of $\varphi^2(\Gamma_{\partial D^2\times D^1})$  align with the two arcs of $\Gamma_{\partial_+W}\cap \varphi^2({\partial D^2\times D^1})$. First,  Legendrian realise the two curves $\partial \varphi^2 ( \partial D^2\times \partial D^1)$ on  $\partial_+W$, and then introduce a corner along each  to obtain a cobordism with piecewise convex boundary.  Finally,  glue $h^2$ using $\varphi^2$.

Contact $3$--handles are easier to attach, as one need not create any corners before gluing. 

When $W$ is embedded in some contact 3-manifold $(M,\xi)$,  one may attach 1-handles to $W$ inside $(M,\xi)$ along any Legendrian arc $l$ properly embedded in $M\setminus W$ with boundary on $\Gamma_{\partial_+W}\subset \partial W$. In this case, the attachment  is the (smoothing of) $\big(W\cup\nu(l),\xi\vert_{W\cup\nu(l)}\big)$, where $\nu(l)$ is a standard neighbourhood of $l$. Up to weak isotopy, this construction depends only on the Legendrian isotopy class of $l$ relative to its endpoints and not on the choice of standard neighbourhood  or the particular representative of the Legendrian isotopy class. 

Attaching a contact $1$--handle preserves tightness:
\begin{lemma} [\cite{Hgluing}, Corollary 3.6(2)] \label{lem:tight1handle} Let $W$ be a contact manifold with a convex boundary and let $W'=W \cup \nu(l)$ be the  manifold formed by attaching a contact $1$-handle.  If $W$ is tight, then so is $W'$.
\end{lemma}

Similarly, given a properly embedded convex disc in $M\setminus W$, one may add a standard or $I$-invariant neighbourhood of the disc to $W$ as a 2-handle attachment.  In order to ensure this is a contact $2$-handle attachment, the disc must have a tight neighbourhood and a Legendrian boundary on $\partial_+W$ with Thurston-Bennequin number $-1$.

As usual, attaching handles to a cobordism changes the boundary by surgery; in this case, contact handle attachment changes $\partial_+W$ by convex surgery, so that the new boundary has a dividing set distinguished up to isotopy by the handle attachment.

\subsection{Bypasses}\label{sec:bypass} 

An isotopy of the convex surface $\Sigma$  in a contact manifold either preserves the dividing set of $\Sigma$ up to isotopy or changes it by a sequence of \textit{bypasses}, each of which corresponds to pushing $\Sigma$ across a particular contact three-ball.  This ball may be characterised in a variety of ways; most familiar is viewing the ball as a neighbourhood of a half an overtwisted disc, but we will use an equivalent definition that is more convenient for our purposes.  

Let $(\Sigma,\Gamma)$ be a convex surface in $(M, \xi)$. An arc $c\subset \Sigma$ is \emph{admissible} if it is transverse to $\Gamma$, $\partial c\in \Gamma$, and the interior of $c$ intersects $\Gamma$ once.  By a $C^\infty$-small convex isotopy of $\Sigma$ the admissible arc $c$ can be made  Legendrian on $\Sigma$; this condition is subsequently assumed without mention.

\begin{figure}[h]
\labellist
\small\hair 2pt
\pinlabel $D$ [l] at 36 75
\pinlabel {$\Gamma_{\Sigma}$} [tr] at 103 77

\endlabellist

\begin{center}
\includegraphics[scale=.95]{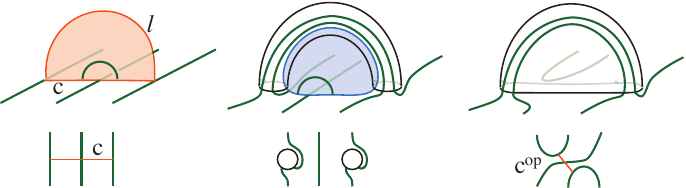}
\caption{ One may isotope $\Sigma$ across the bypass half-disc $D$ cobounded by  $c\cup l$ by successively attaching a contact $1$--handle along $l$ and then a contact $2$--handle along $D$. This decomposition was first noted in Example 5 in \cite{HKM}, and explicit models for the attachments are given in Section 3 of \cite{Oz}.} 
\label{fig:bypassmodel}
\end{center}
\end{figure}

A \emph{bypass disc} $D$ is a convex half-disc with Legendrian and piecewise-smooth  boundary $\partial D=c\cup l$, where 
\begin{enumerate}
\item $D$ intersects $\Sigma$ transversely exactly at $c$; 
\item $D$ has a tight neighbourhood;
\item $\tw_l(\xi,TD)=0$.   
\end{enumerate}
Since $c$ is admissible, $\tw_c(\xi,TD)=-1$.  The fact that $D$ has a tight neighbourhood allows one to draw a dividing curve on $D$ as in Figure \ref{fig:bypassmodel}. 

Suppose that $\Sigma$ is oriented so that $D$ is on its positive side.  We will define a cobordism $W$ that  encloses $\Sigma\cup D$. The operation which replaces the original $\Sigma=\partial_- W$ by $\partial_+ W$ is called \emph{attaching a bypass from the front}.  If $D$ is on the negative side of $\Sigma$, replacing $\Sigma=\partial_-W$ by $\partial_+W$ is called \emph{attaching a bypass from the back}.  

We now construct the cobordism $W$ that encloses $\Sigma\cup D$ for bypass attachment from the front.   
Begin with a standard neighbourhood $\nu(\Sigma)$, where $\Sigma=\partial_-\nu(\Sigma)$.  We also require that  $D\setminus \nu(\Sigma)$ has the same properties as the original $D$, and  we will not distinguish them in the following discussion. Attach the $1$--handle $\nu(l)$ to $\nu(\Sigma)$ and call the resulting contact cobordism $W'$. A further 
 $C^\infty$-small convex isotopy of $\partial_+W'$ ensures that $D\cap \partial_+W'$ is a Legendrian knot. Necessarily,  $D\cap \partial_+W'$ has Thurston-Bennequin number $-1$, so a standard neighbourhood $\nu(D')$ of $D'=D\setminus W'$ is a $2$--handle that we attach to $W'$ to obtain $W$.

The $1$--handle $\nu(l)$ and the $2$--handle $\nu(D')$ are in smoothly  cancelling position, and the $3$--manifold $W$ remains smoothly  isotopic to $\Sigma\times I$.  We call this cobordism a \emph{bypass slice}.  Up to weak isotopy, the contact structure $W$ depends only on  $(\Sigma,\Gamma)$ and the isotopy class of $c$ through admissible arcs. The dividing curve on $\partial_+W$ differs from that of $\partial_-W$ as shown in Figure~\ref{fig:bypassmodel}.
As usual, when  $\Sigma=\partial U$ for a  submanifold $U$, we omit $\nu(\Sigma)$ from our construction and attach the  handles directly to $\partial U$.

Turning a bypass slice upside down exchanges the indices of the $1$-- and $2$--handles, which remain in cancelling position. An upside down bypass slice is thus also a bypass slice, but for attaching the half-disc along the admissible arc $c^{\textrm{op}}$, which is the visible part of the belt sphere of the $1$--handle on $\partial_+W$; see Figure \ref{fig:bypassmodel}. It follows that attaching a bypass from the back along $c^{\textrm{op}}$ is an inverse operation to attaching a bypass from the front along $c$.  See \cite{Oz}, Remark 4.1.

Bypass slices are basic building blocks of contact structures on $\Sigma\times I$.

\begin{theorem}\label{thm:bypass}\cite[Section 3.2.3]{Hgluing}
Any contact structure $\xi$ on $\Sigma\times I$ with convex boundary can be decomposed as a concatenation of bypass slices.  
\end{theorem}

Isotopies are also built up from bypasses in the following sense: 
\begin{theorem}[Colin's Isotopy Discretisation,  \cite{Hgluing} Section 3.2.3 ]\label{thm:isdisc}
Let $\Sigma$ and $\Sigma'$ be convex surfaces in $(M,\xi)$ that are smoothly isotopic. Then there is a sequence of embedded convex surfaces $\Sigma=\Sigma_0,\Sigma_1,\dots, \Sigma_k=\Sigma'$ such that for each $0\le i<k$, the surface $\Sigma_{i+1}$ is obtained from  $\Sigma_{i}$ via a bypass attachment from the front or from the back.  
\end{theorem}

\section{Decompositions of contact manifolds}\label{sec:background2} 
This section discusses convex Heegaard splittings and open book stabilisation.  The statements proven in Section~\ref{sec:obch} are largely a matter of perspective and will be familiar to experts, but Section~\ref{sec:ob} has several technical results that we will rely heavily upon later.

\subsection{Heegaard splittings and open books}\label{sec:obch}

\begin{definition} Let $M$ be a closed orientable $3$-manifold.  The pair $(B, \pi)$ is an (embedded) \emph{open book} if $B$ is an oriented link in $M$ and $\pi:M\setminus B\rightarrow S^1$ is a fibration  which restricts to the normal angular coordinate in a neighbourhood of  $B$.
 \end{definition}

An open book decomposition  $(B,\pi)$ for $M$ defines a Heegaard splitting $\H(B, \pi)$ for $M$ with handlebodies 
\begin{equation*} 
\begin{split}
U& =\overline{\pi^{-1}[0,{1/2}]}, \\
V& =\overline{\pi^{-1}[{1/2},1]}.
\end{split}
\end{equation*}

That the manifolds $U$ and $V$ are indeed handlebodies follows from the fact that the fibre over each point in $S^1$ is an open surface, so $U$ and $V$ are closures of products of a surface with an interval.  Choosing a smooth identification of $\overline{M\setminus \pi^{-1}(1)}$ with $S\times I$ for a model fibre $S$ allows us to easily identify a set of compressing discs for each handlebody.  Writing $S_t$ for the \emph{page} $S\times\{t\}$, let $\varphi_s^t\colon S_s\to S_t$ be the parallel transport fixing $B$ and carrying $S_s$ to $S_t$. A set of properly embedded arcs $\{a_1,\dots,a_k\}$ on $S=S_0$ is an \emph{arc system} if $S \setminus  \{a_1,\dots,a_k\}$ is a (cornered) disc. For any arc system,  the discs 
\begin{equation*} 
\begin{split}
A_{i}& =\bigcup_{t\in[0,{1/2}]} \varphi_0^t(a_i)\\
B_{i}& =\bigcup_{t\in[{1/2},1]} \varphi_0^t(a_i)
\end{split}
\end{equation*}
are compressing discs for $U$ and $V$. Denote the union of $A_{i}$ discs by $\A$ and the union of $B_{i}$ discs by $\B$, respectively.   Here and elsewhere, we will orient the Heegaard surface $\Sigma= \overline{\pi^{-1}(\frac{1}{2}) } \cup \overline{-\pi^{-1}(0)}$ as the boundary of $U$.

Given an open book, Torisu established the existence of a unique positive contact structure $\xi$ such that in the associated Heegaard splitting, $\Sigma$ is convex and $\xi$ restricts tightly to each handlebody \cite{Torisu}.   These properties can be taken to characterise the contact structure supported by an open book.  

\begin{definition}{\cite[Lemma 4.1]{Etnyre}}\label{prop:torisu}  Fix an open book decomposition $(B,\pi)$ for $M$.  Then $\xi$ is \emph{supported by} $(B, \pi)$ if and only if $\H(B, \pi)$ has the following properties:
\begin{enumerate} 
\item $\Sigma$ is convex with dividing curve $\Gamma=B$;
\item $\xi\vert_U$ and $\xi\vert_V$ are both tight. 
\end{enumerate}
\end{definition}

A Heegaard splitting  $\H(B,\pi)$ induced by some open book $(B, \pi)$ is called a \textit{convex Heegaard splitting}, and we show next how to identify convex splittings.  Observe that in $\H(B, \pi)$, the boundary of each disc in $\A \cup \B$ intersects the binding $B$ in 2 points. 
In any handlebody with convex boundary, we call a compression disc whose boundary intersects $\Gamma$  twice essentially a \textit{$\Gamma$ product disc}, and by analogy with arcs on a surface, a set of compression discs which cut a handlebody into a ball is called a \textit{disc system}.  Topologically, a $\Gamma$ product disc system identifies the given handlebody as a product of a surface and an interval;  Torisu shows that the existence of disc systems for both handlebodies characterises convex splittings:

\begin{proposition}\label{prop:critforctct}{\cite[Section 3]{Torisu}}
A Heegaard splitting $(\Sigma, U, V)$  of a contact manifold $(M, \xi)$ is a convex Heegaard splitting if and only if 
\begin{enumerate}
\item $\Sigma$ is convex with dividing curve $\Gamma$;
\item\label{tight} $\xi\vert_U$ and $\xi\vert_V$ are both tight;
\item\label{connected} there exist systems of $\Gamma$ product discs $\A$ for $U$ and $\B$ for $V$.  
\end{enumerate}
\end{proposition}

Cutting a handlebody along a disc system reduces it to a ball with paired discs on the boundary; after cutting along a system of $\Gamma$ product discs, there is a unique way to connect arcs of $\Gamma$ whose endpoints lie on the same disc.  This matches the result of smoothing in the case that the product discs are convex with Legendrian boundary, so we again call the resulting curve $\Gamma$.  In fact, this new $\Gamma$, which now lies on a sphere, is connected:  cutting along an arc system on a page of an open book yields a surface with connected boundary, so cutting $U$ and $V$ along $\A$ and $\B$ yields a connected $\Gamma$.

The  construction of an open book from a convex Heegaard splitting  is compatible with a $1$-parameter family of Heegaard splittings:
\begin{proposition}\label{prop:obisot}
Let $(M,\xi)$ be a contact structure and let $(B_0,\pi_0)$ and $(B_1,\pi_1)$ be open books supporting $\xi$. Then $\H(B_0,\pi_0)$ and $\H(B_1,\pi_1)$ are isotopic via an isotopy keeping the Heegaard surfaces $\Sigma_0$ and $\Sigma_1$ convex if and only if $(B_0,\pi_0)$ and $(B_1,\pi_1)$ are isotopic through a path of open books supporting $\xi$. 
\end{proposition}

\begin{proof} Suppose first that $(B_t, \pi_t)$ is a path of open books supporting $\xi$.  At each point in the path, the induced Heegaard surface $\Sigma_t$ is a union of two pages with dividing set $B_t$, so this defines an isotopy from $\H(B_0,\pi_0)$ to $\H(B_1,\pi_1)$.  Conversely, observe that an isotopy taking $\H(B_0,\pi_0)$ to $\H(B_1,\pi_1)$ carries with it a family of $B_t$ product discs, and hence, defines an open book $(B_t, \pi_t)$ for all $t$.  An isotopy of open book decompositions preserves the supported contact structure, as desired.
\end{proof}

In light of Proposition~\ref{prop:obisot}, we will view open books and Heegaard splittings as objects defined only up to convex isotopy.

A key aim of this paper is to broaden the class of Heegaard splittings that may be effectively used to study contact structures.  To this end, we introduce the following definition: 

\begin{definition}\label{def:tighths} A Heegaard splitting $(\Sigma, U, V)$ of a contact manifold $(M, \xi)$ is  \textit{tight} if $\Sigma$ is convex and $\xi\vert_U$ and $\xi\vert_V$ are both tight.
\end{definition}

Every convex Heegaard splitting is tight, and Proposition~\ref{prop:critforctct} states that a tight splitting is contact if it admits a system of product discs.

\begin{example} Tight splittings are strictly more general than convex splittings.  For example, in a convex splitting, $\Gamma_\Sigma$ divides $\Sigma$ into two connected components, while the figure below shows a tight splitting where $\Sigma\setminus \Gamma_\Sigma$ has four components.

\begin{figure}[h]
\labellist
\small\hair 2pt
\pinlabel $\Gamma_\Sigma$ [l] at 150 90
\pinlabel $\partial B$ [l] at -13 72
\pinlabel $\partial A$ [l] at -13 30

\endlabellist

\begin{center}
\includegraphics[scale=.85]{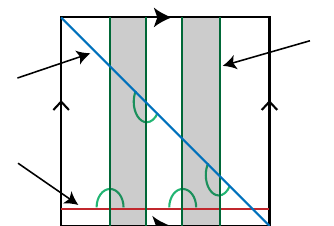}
\caption{  A toroidal Heegaard surface for $S^3$ where $\Gamma_\Sigma$ cuts $\Sigma$ into four components.  The curved green arcs indicate  the arcs of $\Gamma_A$ and $\Gamma_B$. }
\label{fig:disconn}
\end{center}
\end{figure}

Figure~\ref{fig:disconn} shows a convex Heegaard diagram for $S^3$ together with the dividing set on $\Sigma$.  In addition to the line segments, there are curved green arcs indicating which points are connected by the respective dividing sets on the meridional discs.  After cutting either solid torus along the indicated meridian and smoothing the resulting ball, the dividing set is a connected curve.  By Giroux's Criterion (Theorem~\ref{thm:gircrit}), the sphere has a tight neighbourhood and can be filled by a tight ball.  

\end{example}
\subsection{Stabilising open book decompositions}\label{sec:ob} 

Stabilisation is an operation performed on open book decompositions of a $3$--manifold.  First identified by Stallings in the topological context,  stabilisation  comes in two versions that are distinguished as positive and negative \cite{Stallings}.  Positive stabilisation preserves the supported contact structure, up to isotopy, and is the only version considered in this paper.  In this section we establish the equivalence of several perspectives on positive open book stabilisation, focusing on identifying when a change to a Heegaard splitting is in fact a positive stabilisation of the underlying open book.

\begin{figure}[h]
\labellist
\small\hair 2pt
\pinlabel $\gamma$ [l] at 40 25
\endlabellist
\begin{center}
\includegraphics[scale=1.2]{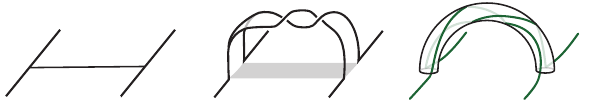}
\caption{ Left: An open book stabilisation is determined by an arc $\gamma$ properly embedded on a page. Centre: Each page of the connect sum  $(B, \pi)\#(H^+, \pi^+)$ is a plumbing between the original page with an annulus at the shaded rectangle.  Right: The genus of the induced Heegaard splitting increases by one.  The dividing set twists as shown to match the twisting of the $1$--handle added to the page.}
\label{fig:stabnew}
\end{center}
\end{figure}

The literature offers several equivalent ways to define positive stabilisation, and we present the one best suited to the later discussion: 
let $(H^+, \pi^+)$ denote the open book for $S^3$ with binding a positive Hopf link $H^+$ in the unit sphere in $\mathbb{C}^2$ and $\pi^+$  the fibration over $S^1$ defined by $(z_1, z_2)\mapsto  \frac{z_1z_2}{|z_1z_2|}$.  Inside $(H^+, \pi^+)$, choose a $3$--ball neighbourhhood of a cocore arc of $(\pi^+)^{-1}(0)$. Given an arbitrary open book $(B,\pi)$, a positive stabilisation is formed by taking the connect sum $(B, \pi)\# (H^+, \pi^+)$ using the designated $3$--ball in $(H^+, \pi^+)$ and a 3--ball in $(B,\pi)$ that is also a neighbourhood of some arc $\gamma$ properly embedded on a page in $M$.    In this case we can arrange that the open book data  match along the gluing $S^2$, and the new page is a Murasugi sum of the two original pages \cite{Gab}, \cite{Etnyre}.  See Figure~\ref{fig:stabnew}. The open book $(B^+,\pi^+)$ obtained thus is a \emph{(positive) stabilisation} of $(B, \pi)$ along $\gamma$.

Now let $\H(B, \pi)$ be the Heegaard splitting induced by an open book decomposition $(B, \pi)$ for $(M,\xi)$.  It is straightforward to check that stabilising $(B, \pi)$ to $(B^+, \pi^+)$ induces a Heegaard splitting stabilisation of $\H(B, \pi)$, as the effect of adding a $2$-dimensional $1$--handle to the page is adding a $3$-dimensional $1$-handle to the handlebody formed by thickening the page. See the right-hand picture in Figure~\ref{fig:stabnew}.  The twist in  the plumbed annulus determines how the dividing set twists around the new handle in the handlebody formed by thickening the page.  

The next definition introduces an alternative process that alters a Heegaard splitting.

\begin{definition}\label{def:posstab}
Suppose that $(\Sigma, U, V)$ is a tight Heegaard splitting of a contact manifold $(M, \xi)$ and let $D\subset V$ be a convex half-disc with Legendrian and piecewise-smooth boundary $\partial D=c\cup l$.   Suppose that $c\subset \Sigma$ and that $l$ is properly embedded in $V$ with $\partial l=-\partial c\subset \Gamma_\Sigma$ and suppose further that 
\[\textit{tw}_l(\xi,TD)=\textit{tw}_c(\xi,TD)=-\frac12.\] 
This data determines a new  Heegaard splitting $(\Sigma', U', V')$  of $M$, where $U'$ is the smoothing of $U\cup \nu(l)$ and  $V'=V\setminus \nu(l)$. We say that the splitting $(\Sigma', U', V')$ is a \emph{positive stabilisation} of $(\Sigma, U, V)$.
\end{definition}

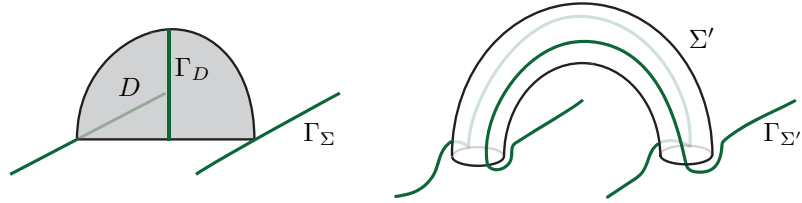
\begin{figure}[h]
\labellist
\small\hair 2pt
\pinlabel $D$ [l] at 60 40
\pinlabel $\Gamma_D$ [l] at 78 46
\pinlabel {$\Gamma_\Sigma$} [tr] at 132 30
\pinlabel {$\Gamma_{\Sigma'}$} [tr] at 280 30
\pinlabel {$\Sigma'$} [tr] at 252 60

\endlabellist

\begin{center}
\includegraphics[scale=1.2]{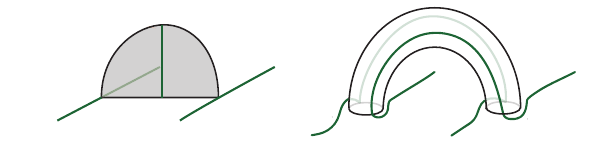}
\caption{ Left: Local model for a half-disc defining a positive stabilisation.  Right: The stabilised Heegaard surface $\Sigma'$ shown with the new dividing set $\Gamma_{\Sigma'}$.}
\label{fig:stab0}
\end{center}
\end{figure}

A local model for $D$ defining a positive stabilisation is shown in Figure~\ref{fig:stab0}.  The fact that $\xi$ is tight on $V$, and thus near $D$, implies that the dividing curve on $D$ is just an arc connecting $c$ and $l$.  Thus $l$ is Legendrian isotopic to $c$. 

The dividing curve on the stabilised surface $\Sigma'$ is dictated by the conditions on $D$.   Observe first that $\textit{tw}=-\frac12$ determines $\Gamma_{\partial \nu(l)}$; then the  edge rounding described in Lemma~\ref{lem:edgerounding} produces $\Gamma_{\Sigma'}$  as shown on the right in Figure~\ref{fig:stab0}.

At first glance, it may seem unfortunate that Definition~\ref{def:posstab} uses a term already in use, but the next lemma redeems this decision:

\begin{lemma}\label{lem:stabHD}
A positive stabilisation of a convex Heegaard splitting $\H(B, \pi)$ is also a convex Heegaard splitting $\H'=\H(B',\pi')$. Furthermore, the open book $(B',\pi')$ is a positive stabilisation of $(B,\pi)$. 
\end{lemma}

\begin{proof}

We first verify that a positive stabilisation of a convex Heegaard splitting is again a convex Heegaard splitting.  The surface $\Sigma'$ is convex by construction, so it remains to verify the second and third conditions listed in   Proposition~\ref{prop:critforctct}.  With respect to the notation introduced in Definition~\ref{def:posstab}, both $D\cap V'$ and the co-core disc of the new handle $\nu(l)$  are product discs for their respective handlebodies.  Cutting along these discs  returns the original Heegaard splitting, so Condition~\ref{connected} is satisfied.  Recall from Section~\ref{sec:handles} that adding $\nu(l)$ to $U$ may be viewed as a contact $1$--handle attachment.   Lemma~\ref{lem:tight1handle} states that contact handle attachment preserves tightness, so $\xi|_{U'}$ is tight.   Since $V'\subset V$, tightness for $V'$ is automatic, and Condition~\ref{tight} is also satisfied.  Thus the stabilised Heegaard splitting is again  a convex Heegaard splitting, as claimed.

To prove the second part of the lemma, view the Legendrian arc $c$ as embedded in the page $\pi^{-1}(\frac{1}{2})$ of $(B,\pi)$.  Perform a positive open book stabilisation along $c$ and call the new open book $(B',\pi')$.  The right-hand picture in Figure~\ref{fig:stabnew} shows the local change to the Heegaard surface for $\H(B',\pi')$, and we see that this is convexly isotopic to the positively stabilised Heegaard surface shown in Figure~\ref{fig:stab0}.
\end{proof}

\begin{remark}  Positive stabilisation replaces one Heegaard splitting of a contact manifold with another Heegaard splitting of the same contact manifold.  There is a natural topological operation that one might be tempted to call negative stabilisation: increase the genus of the Heegaard splitting by one,  but twist the dividing curve in the opposite direction around the new handle.  In light of Lemma~\ref{lem:stabHD}, one may see that this corresponds to the Heegaard splitting induced by plumbing each page of the original open book with an annulus that twists in the opposite direction.  Such an open book is known as a \textit{negative stabilisation} of the original open book and does not support the same contact structure.  
\end{remark}

Open book decompositions of a fixed manifold are partitioned into positive stabilisation classes.  In particular, the property that two open book admit a common positive stabilisation is an equivalence relation.  This is  easily seen by noting that the connect sum stabilising an open book is taken by removing an arbitrarily small $3$--ball neighbourhood of an arc on a page.  Consider two sequences of positive stabilisations of a fixed open book as a sequence of arcs on ordered pages.  Then one may construct a common positive stabilisation by using representatives of all the arcs on any set of ordered pages that restricts to the correct order on each of the two subsequences.

 We conclude this section with a result about bypass attachment which will be essential  in the later sections:

\begin{proposition}\label{prop:bypass}  Suppose $\H=(\Sigma, U, V)$ is a convex Heegaard splitting of the a contact manifold $(M,\xi)$ and suppose that $\H'=(\Sigma', U', V')$ is a new convex splitting, where $\Sigma'$ is obtained from $\Sigma$ by a bypass attachment along the bypass half-disc $D$ in $V$. If $D$ is disjoint from a system of product discs $(\A, \B)$ for $\H$, then  $\H$ and $\H'$ admit a common positive stabilisation.  

An analogous statement holds if the bypass half-disc is in $U$. 
\end{proposition}

\begin{proof}[Proof of Proposition~\ref{prop:bypass}] As before, label the Legendrian arcs of $\partial D$ as $c\subset \Sigma$ and $l \subset V$.  Let $\H''$ be the convex Heegaard splitting built by transferring a contact $1$--handle neighbourhood of $l$ from $V$ to $U$.  We will show that $\H''$ is a positive stabilisation of $\H$ and $\H'$.    

We begin by examining the bypass attachment arc $c\subset \Sigma$.  Since $c$ crosses $\Gamma_\Sigma$, it does not satisfy the hypotheses of Lemma~\ref{lem:stabHD} that ensure positive stabilisation, but we claim that the hypothesis $D\cap\B=\emptyset$ implies the existence of  an alternative disc $D'$ with boundary $c'\cup l$ which does satisfy the conditions of Lemma~\ref{lem:stabHD}. Such a $D'$ suffices to prove that $\H''$ is a positive stabilisation of $\H$.

In the following we will concentrate on the existence of $D'$. The contact handlebody $V$ is obtained from a ball $B_V$ via contact 1-handle attachments with co-cores $B_i$, as indicated on the left-hand side of Figure~\ref{fig:bypstab2}.  Since $D$ lies in the complement of the discs $\B$,   attaching  a neighbourhood of $D$ to $U$  is a trivial bypass in $B_V$; again, see Figure~\ref{fig:bypstab2}.  The alternative disc $D'$ can be found in $B_V\subset V$ as shown on the right-hand side of Figure~\ref{fig:bypstab2}. Explicitly,  construct $D'$ from $D$ by sliding  $D$ across the co-cores of the 1-handles that meet the indicated subarc $\gamma$ of $\Gamma_\Sigma$.

\begin{figure}[h]
\labellist
\small\hair 2pt
\pinlabel $\textit{l}$ [l] at 85 83
\pinlabel $\textit{l}$ [l] at 245 83
\pinlabel $D$ [l] at 90 65
\pinlabel $c$ [l] at 90 40
\pinlabel $D'$ [l] at 250 65
\pinlabel {$c'$} [tr] at 245 40
\pinlabel {$\gamma$} [tr] at 70 30

\endlabellist

\begin{center}
\includegraphics[scale=1.2]{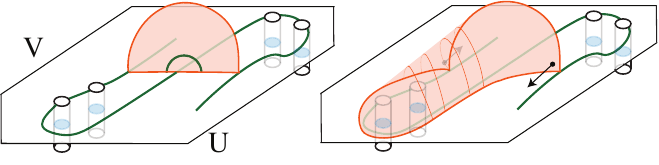} 
\caption{  Left: Here, $V$ is shown with indicative product discs.  Observe that after cutting along all the product discs and smoothing, the new dividing curve $\Gamma_{S^2}$ is connected, so the bypass is trivial in the resulting tight ball.  Right: Isotoping $c\cap \overline{-\pi^{-1}(\frac12)}$ across the product discs yields a curve $c'$ which is disjoint from $\Gamma_\Sigma$ and cobounds a new half-disc with $l$.} 
\label{fig:bypstab2}
\end{center}
\end{figure}
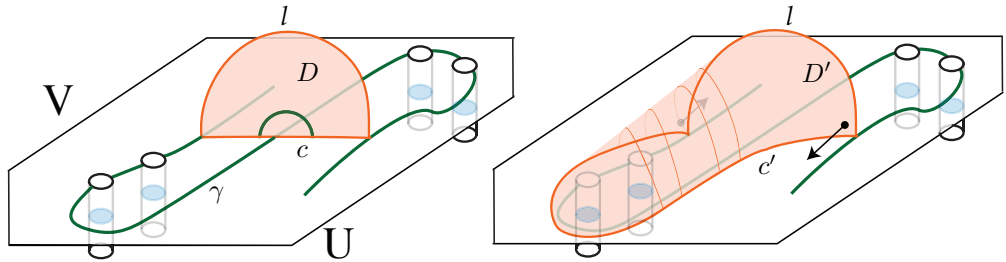

 We confirm that $D'$ defines a positive stabilisation.  Since $D$ is a bypass half-disc, the twisting of $\xi$ along $l$ relative to $D$ is $0$.  The half-discs coincide near the endpoint where $c=c'$,   but the orientation of $D'$ is opposite that of $D$ at the other shared endpoint.   This implies that the relative twisting of $\xi$ and $D'$ is non-zero, but it cannot exceed $-\frac{1}{2}$ because the interiors of $D$ and $D'$ may be isotoped to be disjoint, by construction.  Thus $\tw_l (\xi, D')=-\frac{1}{2}$.

The other boundary component of $D'$ is the curve $c'$ properly embedded in  $\overline{-\pi^{-1}(0)}\subset \Sigma$.  Since $\Sigma$ is convex and $c'\cap \Gamma_{\Sigma}=\partial c'$, it follows that $\tw_{c'} (\xi, D')=-\frac{1}{2}$ as well.  This establishes that $\H''$ is a positive stabilisation of $\H$.

It remains to show that $\H''$ is a positive stabilisation of $\H'$, as well.  Recall the duality of bypass attachment from the front and back discussed in  Section~\ref{sec:bypass} : $U$ can be obtained from $U'$ by attaching a bypass $D^\textrm{op}$ along $c^\textrm{op}$.  As $D$ was disjoint from $\A\cup\B$, the opposite bypass disc $D^\textrm{op}$ is also disjoint from $\A\cup\B$.  Attaching the $1$--handle corresponding to the bypass $D^\textrm{op}$ once again produces $\H''$, and a parallel argument establishes that this is a positive stabilisation.  
\end{proof}

\section{Existence}\label{sec:ex}
Here we use the language of contact handle decompositions to present a proof for the existence of an open book compatible with any contact structure:

\begin{theorem}\label{thm:ex2}
Let $(\Sigma, U, V)$ be a Heegaard splitting of a contact $3$--manifold $(M, \xi)$.  There exists a sequence of Heegaard splitting stabilisations such that the resulting Heegaard splitting  $(\Sigma', U', V')$ is isotopic to a convex splitting for $(M, \xi)$.
\end{theorem}

As convex Heegaard splittings and open books are equivalent, the following result is an immediate consequence of Theorem~\ref{thm:ex2}.

\begin{corollary}[c.f.\ Theorem \ref{thm:ex}]\label{cor:ex} Any contact $3$--manifold  admits a supporting open book decomposition. 
 \end{corollary}

\begin{proof}[Proof of Theorem~\ref{thm:ex2}]
Take any (smooth) Heegaard splitting $(\Sigma, U, V)$ of $M$ and choose bouquets of circles $K_U\cong \vee_{i=1}^g S_i^1$ inside $U$ and $K_V\cong \vee_{j=1}^g S_j^1$ inside $V$ such that each handlebody deformation retracts onto its respective bouquet. 
Legendrian realise each of $K_U$ and $K_V$. Retaining the same label, take  standard  neighbourhoods $\nu(K_U)$ and $\nu(K_V)$. Each of these is  a tight handlebody with a convex boundary, but their union doesn't exhaust $(M,\xi)$.

 The closure of the complement $(M\setminus \nu(K_U \cup  K_V),\xi\vert_{M\setminus \nu(K_U\cup  K_V)})$ 
 is a contact manifold diffeomorphic to $\Sigma\times I$ with convex boundary. By Theorem \ref{thm:bypass}, it decomposes into, say, $k$ bypass slices. 
As described in Section~\ref{sec:bypass}, each bypass slice in turn decomposes into a contact 1-handle $(h_i^1,\zeta^1_i)$ and a contact 2-handle $(h_i^2,\zeta^2_i)$.  Thus $(\Sigma\times I,\xi)$ is weakly contact isotopic to the smoothed \[(\nu(\Sigma_U),\xi_U)\cup \bigcup_{i=1}^{k} \left((h_i^1,\zeta^1_i)\cup(h_i^2,\zeta^2_i)\right),\]
where $(\nu(\Sigma_U),\xi_U)$ is an $I$-invariant half-neighbourhood of $\Sigma_U:=\partial \nu(K_U)$. 
Moreover, for dimension reasons one can assume that the attaching region of any 1-handle $h_j^1$ or 2-handle   $h_j^2$ is disjoint from $\cup_{i<j} (\partial h_i^2,\zeta^2_i)$. This means that the bypass attaching arc $c_i$ is also disjoint from $\cup_{i<j} (\partial h_i^2,\zeta^2_i)$.

Consider \[U'=\nu(K_U)\cup \bigcup_{i=1}^{k} h_i^1 \]
and set $V'=M\setminus U'$. We claim that $M=U'\cup V'$ is a convex Heegaard splitting of $(M,\xi)$. Indeed, $U'$ is  tight and  decomposes along  product discs to a ball with a connected dividing set, as both these properties persist under $1$--handle addition. To make a similar claim for $V'$, we turn the picture upside down.  Given
\[V'=\nu(K_V)\cup \bigcup_{i=1}^{k} h_i^2, \]
 we view the $h_i^2$ as $1$--handles  attached to $\partial \nu(K_V)$ in the opposite order. Thus $V'$ is also tight and admits a system of product discs, as desired.

Examining  the construction with ``smooth glasses'', observe that the Heegard splitting $\big(\Sigma_U, \nu(K_U), M\setminus  \nu(K_U)\big)$
 is smoothly isotopic to the original $(\Sigma, U,V)$.  If we build $U'$ by sequential $1$--handle additions,  each step is a Heegaard splitting stabilisation of the previous one.  This is easily seen, since the  $2$--handle $h_i^2$ cancels the $1$--handle $h_i^1$ in the smooth category.   Thus $(\Sigma', U', V')$ is indeed a Heegaard splitting stabilisation of $\H$.
\end{proof}

\section{Tight Heegaard splittings}\label{sec:tighths} 
In this section we show that even non-convex Heegaard splittings of a tight contact manifold can be used to produce  open books.  A tight Heegaard splitting (Definition~\ref{def:tighths}), together with certain systems of compressing discs, allows us to construct a new convex Heegaard splitting, and thus, an open book.   This process is called \emph{refinement.}    We show that the positive stabilisation class of the resulting open book depends only on the original tight Heegaard splitting of $(M, \xi)$. (See Theorem~\ref{thm:indep}).

\subsection{Refinement: Special Case} 

We first consider the special case where  $\Sigma$ is decorated with Legendrian attaching curves for a pair of disc systems for the two handlebodies.  More precisely,  let $\H=(\Sigma, U, V)$ be a tight Heegaard splitting of $(M, \xi)$.  Let $\A=\{A_1,\dots,A_g\}$ be a disc system for $U$ such that for each $i$, $\partial A_i$  
is Legendrian, $A_i$ is convex, and the intersection $A_i\cap\Sigma$ is standard.   Let  $\B=\{B_1,\dots,B_g\}$ be a  disc system for $V$ with the same properties.  We further assume that the multicurves  $\partial\A$, $\partial\B$, and  $\Gamma_\Sigma$ are in general position on $\Sigma$.  We call the pair $(\A, \B)$ a  \emph{convex compressing disc system} for $\H$.  Note that on an arbitrary convex Heegaard surface, it need not be possible to simultaneously Legendrian realise $\partial \A \cup \partial \B$; this more general situation is addressed in Section~\ref{sec:gen}.

Given a tight Heegaard splitting $\H$ with a convex compressing disc system $(\A,\B)$, we will show how to construct a convex Heegaard splitting $\widetilde{\H}(\A,\B)=(\widetilde{\Sigma},\widetilde{U}, \widetilde{V})$ 
of $(M,\xi)$, and thus, an open book for $(M,\xi)$.  The convex Heegaard splitting $\widetilde{\H}(\A,\B)$ is called the \emph{contact refinement} of $\H$ via $(\A,\B)$.  

Roughly speaking, we construct the refinement  corresponding to $(\A,\B)$ by tunnelling along a spine of $A_i\setminus \Gamma_{A_i}$ and $B_j\setminus \Gamma_{B_j}$ \footnote{For example, on the characteristic foliation one can take the union of the (usual) graphs formed by the positive and negative singularities, separately. } and transferring each of the excavated $1$-handles to the opposite handlebody.  This has the effect of breaking the $A_i$ and $B_j$  into  collections of product discs in a new convex Heegaard splitting.  Let us now describe this more precisely.

Suppose that $(\A, \B)$ is a convex compressing disc system for the splitting $\H$.  Note that $\Gamma_{\Sigma}$ cannot be empty;  as $\partial A_i$ and $\partial B_j$ are Legendrian, both $\Gamma_{A_i}$ and $\Gamma_{B_j}$ are also non-empty.  

\begin{lemma}\label{lem:nondisj} Each component of  $\Gamma_{A_i}$ intersects $\partial A_i$ and each component of  $\Gamma_{B_j}$ intersects $\partial B_i$.  Each component of $\Gamma_\Sigma$ intersects both $\partial \A$ and $\partial \B$.  Each $\partial A_i$ and each $\partial B_j$ intersects $\Gamma_\Sigma$.
\end{lemma}

\begin{proof} If  any component $\gamma$ of $\Gamma_{A_i}$ or $\Gamma_\Sigma$ were disjoint from $\partial \A$, then $\gamma$  would persist after cutting along $\A$ and smoothing.  This would yield a disconnected dividing set on a tight $3$--ball, contradicting Theorem~\ref{thm:gircrit}.  Similarly,  any component of $\Gamma_{B_j}$ or $\Gamma_\Sigma$  disjoint from $\partial \B$ would persist after cutting and smoothing and again lead to a contradiction. 

For the final claim, observe that any component $ A_i$  or $B_j$ disjoint from $\Gamma_\Sigma$ is necessarily an overtwisted disc.
\end{proof}

On each disc  $A_i$, endpoints of $\Gamma_{A_i}$ and $\Gamma_{\Sigma}$ alternate along ${\partial A_i}$.  Let $X_i^\A$ be a collection of properly embedded arcs in $A_i\setminus \Gamma_{A_i}$ with disjoint interiors and endpoints on  ${\partial A}_i \cap \Gamma_{\Sigma}$ such that $X_i^\A$ cuts $A_i$ into subdiscs each containing a single component of $\Gamma_{A_i}$. Using the Legendrian Realisation Principle,  perform a  $C^\infty$-small convex isotopy of $A_i$ relative to $\partial A_i$ to ensure that the arcs $X_i^\A$ are Legendrian on $A_i$.  Let $\nu(\X^\A)$ be a standard contact neighbourhood of the union of these arcs and define $\overline{V}$ to be the smoothing of $V\cup \nu(\X^\A)$.  Let $\overline{U}$ be $M\setminus \overline{V}$. 

Taking the union of $V$ and $\nu(\X^\A)$ requires standard intersections. To achieve this, one must first isotope $\Sigma$ to make the intersections ${\partial }\nu(\X^\A)\cap\Sigma$ Legendrian.  With this achieved, $\nu(\X^\A)$ must be isotoped relative to $\Sigma$ to ensure   that the intersection between $\partial \nu(\X^\A)$ and $\Sigma$ is standard.   Then $V\cup\nu(\X^\A)$ can be smoothed to have a convex boundary.   Since this is both possible to do and painstaking to describe, we may suppress such details in the following. 

Continuing, let $\X^\B$ be an analogous collection of arcs cutting $\B$ into subdiscs each containing one component of $\Gamma_{\B}$.
Repeat this process to produce new handlebodies $\widetilde{U}=\overline{U} \cup \nu(\X^\B)$ and $\widetilde{V}=M\setminus \widetilde{U}$.  The resulting Heegaard splitting $\widetilde{\H}=(\widetilde{\Sigma}, \widetilde{U}, \widetilde{V})$ is the \emph{refinement of $\H$ via $\A,\B,$ and $\X$}. 

 \begin{figure}[h]
 \labellist
\small\hair 2pt
\pinlabel {$A_i$} [tr] at 15 110
\pinlabel {$\nu(\Gamma_\Sigma \cap \partial A_i)$} [tr] at 165 38
\endlabellist
\begin{center}
\includegraphics[scale=.95]{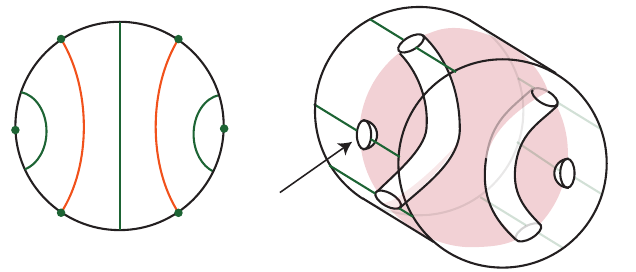}
\caption{ Left: A compressing disc $A_i$ with its green dividing set $\Gamma_{A_i}$ and the orange arcs of $\X$ that separate components of $\Gamma_{A_i}$.    Right: When $A_i$ is part of a convex compressing disc system, remove a standard neighbourhood of the $\X$ arcs  to refine the splitting.   The  compressing discs for the new splitting are the remaining shaded regions.  As shown, we may also remove neighbourhoods of  points of $\Gamma_\Sigma\cap \partial A_i$ that lie in bigons cut out by $\Gamma_{A_i}$.} 
\label{fig:Xfig}
\end{center}
\end{figure}

Later, we will find it convenient to remove a neighbourhood of every point of $\partial \A_i \cap \Gamma_\Sigma$ from $U$.  When such points lie in bigons cut out by $\Gamma_{A_i}$, as opposed to being endpoints of arcs of $\X$, then this changes $\widetilde{\Sigma}$ only by a convex isotopy.

Refinement produces  convex Heegaard splittings:

\begin{proposition}\label{prop:virtcont} The refinement  of $\H$ via $\A,\B,$ and $\X$ is a convex Heegaard splitting of $(M, \xi)$.
\end{proposition}

\begin{proof}  The Heegaard surface $\widetilde{\Sigma}$ is convex by construction, and we will show that each of the new handlebodies $\widetilde{U}$ and $\widetilde{V}$ is  tight and admits a system of convex product discs. By Proposition~\ref{prop:critforctct}, this suffices to show  $(\widetilde{\Sigma}, \widetilde{U}, \widetilde{V})$ is a convex Heegaard splitting.

Recall that $\overline{V}$ denotes the intermediate handlebody $V\cup \nu(\X^\A)$ and $\overline{U}:=M\setminus \overline{V}$.  By hypothesis, smoothing $U\setminus \A$ yields a tight ball, and the handlebody $\overline{U}$ is obtained from this ball via contact contact 1-handle attachments whose co-cores are the product discs $\A\setminus \nu(\X^\A)$.  The final handlebody $\widetilde{U}$ is obtained by further attaching the contact $1$--handles $\nu(\X^\B)$; each of these again admits a co-core product disc intersecting $\Gamma_{\widetilde{\Sigma}}$ twice.  It follows that cutting along the new product discs again yields a tight ball, and hence one with a connected dividing set.  Since $\widetilde{U}$ is reconstructed from a tight ball by contact $1$-handle attachments, it follows from Lemma~\ref{lem:tight1handle} that $\widetilde{U}$ is tight.  An identical argument applies to $\widetilde{V}$. 
\end{proof}

In fact,  the choice of arcs $\X= \X^\A\cup\X^\B$ does not affect the convex splitting.

\begin{lemma}\label{lem:indx} Let $\X$ and $\X'$ be two sets of arcs separating components of $\Gamma_{\A \cup \B}$.   Then the refinement of $\H$ via $\A, \B$, and $\X$ and the refinement of $\H$ via $\A, \B$, and $\X'$ are convexly isotopic.

\end{lemma}

\begin{proof}
On each disc, the $\X$ and $\X'$ arcs are related by a sequence of arc slides. To see this, one may easily verify that arc slides suffice to transform an arbitrary set of $\X$ curves on $A_i$ into one where all the arcs share a single endpoint, and any two such configurations are also related by arc slides. Given this, it suffices to show that a single arc slide preserves $\Gamma_{\widetilde{\Sigma}}$ up to convex isotopy.

 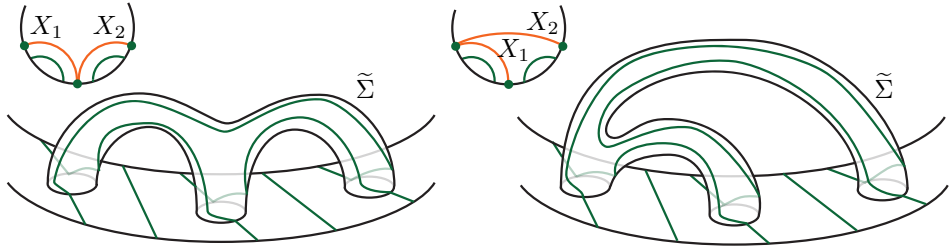
\begin{figure}[h]
 \labellist
\small\hair 2pt
\pinlabel {$X_1$} [tr] at 32 110
\pinlabel {$X_2$} [tr] at 60 110
\pinlabel {$X_1$} [tr] at 239 100
\pinlabel {$X_2$} [tr] at 252 112
\pinlabel {$\widetilde{\Sigma}$} [tr] at 170 85
\pinlabel {$\widetilde{\Sigma}$} [tr] at 400 85
\endlabellist

\begin{center}
\includegraphics[scale=.85]{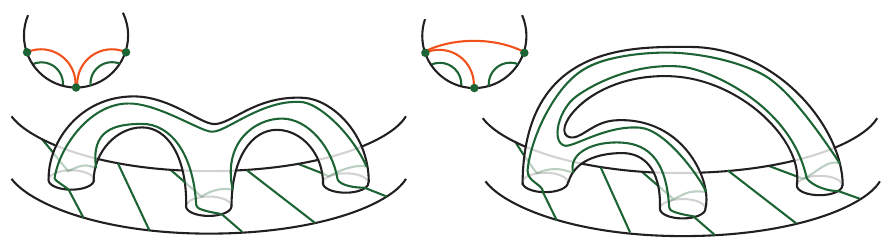}
\caption{ Inset: Two configurations for $\X$  related by a handle slide of $X_2$ over $X_1$.  Main figures: For each of the inset configurations, the smoothed Heegaard surface $\widetilde{\Sigma}$ is shown with its dividing set. } 
\label{fig:xindep}
\end{center}
\end{figure}
Let $X_1$ and $X_2$ be two $\X$ arcs that share an endpoint. We examine a standard neighbourhood $\nu(X_1\cup X_2)$.  
Each arc $X_i$  is disjoint from $\Gamma_\A$, so $\tw_{X_i}(\xi,T\A)=0$.  Thus the dividing set on $\partial \nu(X_1\cup X_2)$ is isotopic relative to its boundary to $\nu(X_1\cup X_2)\cap \A$.
  After smoothing, the curves of $\Gamma_{\partial \nu(X_1\cup X_2)}$ connect to those of $\Gamma_\Sigma$ following Lemma \ref{lem:edgerounding}. Figure~\ref{fig:xindep} shows the local model for the smoothed $\widetilde{\Sigma}$  for two $\X$ configurations related by a single arc slide. It is easy to verify that their dividing sets are isotopic, so the refinements are convexly isotopic.
\end{proof}

Since any refinement of $\H$ is a convex Heegaard splitting, each set of convex compressing discs for $\H$ induces an open book decomposition of the original contact manifold.  Before investigating how this open book depends on the ancillary data, we turn our attention to the case of topological Heegaard splittings where the boundaries of the compressing discs are not Legendrian.

\subsection{Refinement: General Case}\label{sec:gen} 
In the previous section, we began with a tight splitting and a pair of disc systems whose boundaries were simultaneously Legendrian on $\Sigma$. We now relax this condition and consider 
convex Heegaard surfaces where the curves $\partial \A$ are Legendrian realisable and, separately, the curves $\partial \B$ are Legendrian realisable. We want to define the refinement in this case, as well.   In order to do so, we consider a parallel copy $\Sigma'$ of $\Sigma$, chosen so that $\partial \A$ is Legendrian on $\Sigma$ and $\partial \B$ is Legendrian on $\Sigma'$.

\begin{definition} A \emph{triple decomposition}  $\underline{\H}=(N, U, V')$ \emph{underlying} the Heegaard splitting $\H=(\Sigma, U, V)$  is a decomposition of $M$ into pieces $U$, $N$ and $V'$, where $(N,\xi\vert_N)$ is weakly isotopic to an $I$-invariant half-neighbourhood of $\Sigma$. The isotopy is relative to $\Sigma=\Sigma\times\{0\}$ and a neighbourhood of $\Gamma_{\Sigma}\times I$. Set $\Sigma'=\partial_+ N=-\partial V'$.
\end{definition}

\begin{definition} Suppose that $\H$ is a tight Heegaard splitting with convex splitting surface $\Sigma$.  A \emph{convex compressing disc system} $(\A,\B)$ for 
$\H$ is a set of convex compressing discs $\A=\{A_1,\dots,A_g\}$, each with Legendrian boundary on $\Sigma$, and similarly, a set of convex compressing discs $\B=\{B_1,\dots, B_g\}$ with Legendrian boundary on $\Sigma'$ for some underlying triple decomposition $\underline{\H}$.  
\end{definition}

In fact, convex compressing disc systems are common. Given a convex  Heegaard surface, any set of smooth compressing discs for the two handlebodies gives rise to a system of convex compressing discs.  Let  $\A^s=\{A_1^s,\dots, A_g^s\}$ and $\B^s=\{\B_1^s,\dots, B_g^s\}$ be systems of smooth compressing discs for $U$ and $V$, respectively, and assume that $\partial \A^s$ and $\partial \B^s$ are both in general position with $\Gamma_\Sigma$.  Since any system of compressing discs is non-separating, $\partial \A^s$ is non-isolating on $\Sigma$ and may be Legendrian realised; this requires a convex isotopy of the original Heegaard surface, but we again call the resulting surface $\Sigma$.  As a next step isotope $\A^s$ relative to its boundary to be convex and call it $\A$.
  Take a half neighbourhood $\nu(\Sigma)$ of $\Sigma$ and consider the curves $\partial \B^s$ on $\Sigma\times \{1\}$. Slightly abusing  notation, we still denote $\B^{s}\cap (V\setminus \nu(\Sigma))$  by $\B^s$.  Legendrian realise $\partial \B^s$ by isotoping $\Sigma\times \{1\}$ and call the result $\Sigma'$.
 This isotopy can be achieved in an arbitrarily small neighbourhood of $\Sigma\times \{1\}$; $\Sigma'\cap \Sigma=\emptyset$; and we may assume that a neighbourhood of $\Gamma_{\Sigma\times \{1\}}$ is fixed throughout the process.  As a last step, isotope $\B^s$ to be convex and call it $\B$.

We use the isotopy between $N$ and $\nu(\Sigma)$ to ``project'' the curves $\partial \B$ on $\Sigma\times \{1\}$ to curves on  $\Sigma=\Sigma\times\{0\}$ that  need not be Legendrian.  Denote these projected curves on $\Sigma$ by $\partial \B$ as well and extend $\B$ across $\Sigma \times I$ accordingly. The projection naturally identifies the intersection points $\partial \B\cup \Gamma_{\Sigma'}$ on $\Sigma'$ with the intersection points $\partial \B\cup \Gamma_{\Sigma}$ on $\Sigma$,but we are free to alter alter the projection to change  $\partial \A \cap \partial \B$ arbitrarily without changing anything in the upcoming construction. For simplicity, we always assume that the triple $(\partial \A,\partial\B, \Gamma_\Sigma)$ is in general position. 

To recover the special case introduced first, take $N$ to be an $I$-invariant half neighbourhood of $\Sigma$ instead of merely being weakly isotopic to such a neighbourhood.

\subsubsection{Refinement} The refinement in the general case is similar to the special case above, with modifications only to account for the role of  $N$. 

As above, we choose $\X^\A$ and $\X^\B$ and Legendrian realise them on $\A$ and $\B$, respectively. Fix the Legendrian arcs of $\X^\A$, but extend $\X^\B$ through 
$N$ by the curves $(\Gamma\cap\partial \B)\times I$. Since a neighbourhood of $\Gamma\times I$ is fixed by the isotopy bringing $N$ to $\Sigma\times I$, it follows that the arcs $\partial \B \times I$ are automatically Legendrian; we may ensure that the extended arcs are smooth by choosing  $\X^\B$ to be ``straight'' near $\Sigma'$.  Denote the extended arcs by $\X^\B$ so that the refining process proceeds verbatim: for the intermediate splitting, $\overline{V}$ is the smoothed $V'\cup \nu(\X^\A)$ and $\overline{U}=M\setminus \overline{V}$. Finally, $\widetilde{U}=\overline{U}\cup\nu(\X^\B)$ and $\widetilde{V}=M\setminus\widetilde{U}$. The obtained Heegaard splitting  $\widetilde{\H}$ is the \emph{refinement of $\H$ via $\A, \B$ and $\X$}.  
 
 \begin{proposition}\label{prop:vertisctct}  The refinement of $\H$ via $\A, \B$ and $\X$ is a convex Heegaard splitting of $(M, \xi)$.  Furthermore, suppose $\X$ and $\X'$ are two sets of arcs separating components of $\Gamma_{\A \cup \B}$.  Then the refinement of $\H$ via $\A, \B$, and $\X$ and the refinement of $\H$ via $\A, \B$, and $\X'$ are convexly isotopic.
 \end{proposition}
 
 In light of Proposition~\ref{prop:vertisctct}, we drop the reference to $\X$ when describing a refinement.

 \begin{proof}  The proof of Proposition \ref{prop:virtcont} applies directly to show that $U$ and $V'$ are tight handlebodies each cut into a ball by a system of product discs.  To see that the criteria of Proposition~\ref{prop:critforctct} are met, we set $\Sigma=\partial U$ to be the convex splitting surface.  Extending $V'$ across $N$ to $\Sigma$ naturally extends the existing product discs by a collar neighbourhood in $N$, but they remain product discs.  To show that this extended $V'$ is tight, it suffices to note that $N$ is weakly isotopic to a product neighbourhood of a convex surface.

The proof of Lemma \ref{lem:indx} applies verbatim to the general setting.
   \end{proof}

As often in convex surface theory, we see that the smooth object is sufficiently determined by combinatorial input.  In this case, observe that the open book constructed via refinement depends only on the combinatorics of the dividing sets on the two-complexes $\Sigma \cup \A$ and $\Sigma'\cup \B$.  Of note, the intersections between $\partial \A$ and $\partial \B$ are immaterial.  We regard maintaining two surfaces  $\Sigma$ and $\Sigma'$ as a technicality, rather than an essential feature, and when it is unlikely to cause confusion, we may simply write $\Sigma$ both for $\Sigma$ and $\Sigma'$ and $V$ instead of $V'$.

\begin{lemma}\label{lem:refstab} Suppose that  $(\A, \B)$ is a convex compressing disc system for the tight Heegaard splitting $\H$.  If  $\H'$ is a positive stabilisation along a disc $D$ in the complement of $(\A, \B)$, then the refinements of $\H$ and $\H'$ via $\A$ and $\B$ differ by a positive stabilisation. 
\end{lemma}

Lemma~\ref{lem:stabHD} then implies that the open books associated to the two refinements differ by a positive open book stabilisation.

\begin{proof}
This lemma follows from observing that, because $D$ and $\A\cup \B$ are mutually disjoint, one may perform the associated Heegaard splitting stabilisations in any order. 
\end{proof}
 
Combining Lemma~\ref{lem:refstab} with Proposition~\ref{prop:bypass} then yields the following:
 
 \begin{corollary}\label{cor:refstab} Suppose $\H=(\Sigma, U, V)$ and $\H'=(\Sigma', U', V')$ are tight Heegaard splittings related by single bypass attachment to the front or back of $\Sigma$.  If $(\A, \B)$ is a convex compressing disc system disjoint from the bypass half-disc $D$, then the refinements of $\H$ and $\H'$ via $\A$ and $\B$  admit a common positive stabilisation. 
 \end{corollary}

\begin{remark}There is an alternative to constructing a triple as described above.  Given  a tight Heegaard splitting $(\Sigma, U,V)$ for $(M, \xi)$, suppose now that $\A$ and $\B$ are disc systems for the handlebodies satisfying only the requirement that $\partial \A, \partial \B,$ and $\Gamma_\Sigma$ are in general position on $\Sigma$.  If the graph $\partial \A \cup \partial \B$ is non-isolating, then an application of the Legendrian Realisation Principle ensures a convex isotopy of $\Sigma$ that renders the graph Legendrian.  After a further isotopy of the discs relative to their boundaries, $(\A, \B)$ may be assumed a convex compressing disc system.  In the case that the original $\A \cup \B$ is isolating, we claim that it is always possible to perform topological finger moves on $\partial \A$ and  $\partial \B$ to produce a non-isolating graph $\partial\A' \cup \partial\B'$.  One must then show that, up to positive stabilisation, the open book constructed via refinement is independent of the choice of finger moves made at this initial step.  Although this is possible, the flavour of argument is rather different from the rest of the paper, so we have chosen to restrict our discussion to the approach above.
\end{remark}

\subsection{Invariance of the refinement}\label{sec:inv}

Above, we established that refinement promotes a tight Heegaard splitting to a convex Heegaard splitting.  A priori, the latter depends on the choice of convex compressing system, but in this section we show that different choices  preserve the positive stabilisation class of the resulting open book.

\begin{theorem}\label{thm:indep}

Let $(\A,\B)$ and $(\A',\B')$ be two convex compressing disc systems for the tight Heegaard splitting $\H$ of $(M,\xi)$. Then the refinements
$\H(\A,\B)$ and $\H(\A',\B')$ admit a common positive stabilisation. 
\end{theorem}

In order to prove Theorem~\ref{thm:indep}, we will introduce some moves relating distinct convex compressing disc systems.  Each move discretely changes the combinatorics of $(\partial \A,\partial\B,\Gamma_\Sigma)$ and the dividing curves on $\A$ and $\B$.  After introducing each move, we will show that the refinements associated to the two systems are related by a sequence of positive stabilisations.  Finally, we conclude in Proposition~\ref{prop:elmoves} that the  moves considered here suffice to relate any pair of convex compressing disc systems for a fixed tight splitting.

Throughout, let $\H=(\Sigma, U, V)$ be a tight Heegaard splitting for $(M,\xi)$ and let $(N, U, V')$ be an underlying triple decomposition. 
The first two moves occur in an $I$-invariant neighbourhood of $\Sigma$ (or $\Sigma$').

\subsubsection*{[T]: Triple point move} 

\begin{figure}[h]
\labellist
\small\hair 2pt
\pinlabel $An$ [l] at 18 85
\pinlabel $\partial \B$ at 148 70
\pinlabel $\Gamma_\Sigma$ at 140 23
\endlabellist

\begin{center}
\includegraphics[scale=1]{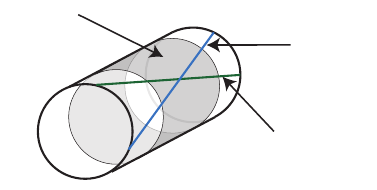}
\caption{ Topological model for the triple point [T] move on a convex compressing disc system.  } \label{fig:moveanst}
\end{center}
\end{figure}

Let $(\A,\B)$ be a convex compressing disc system and let $An \subset \Sigma$ be an annulus with Legendrian boundary $\partial A_1\cup \alpha$. 
Suppose  that  $(\partial \B \cup \Gamma)\cap An$ consists of  arcs connecting $\partial A_1$ to $\alpha$.  These arcs are required to be parallel, with the exception of a single component of  $\partial \B \cap An$ and a single component of $ \Gamma_\Sigma \cap {An}$ which cross once. 
Define $A_1'$ to be a convex surface properly embedded in $U$ that is convex isotopic, relative to $\alpha$, to the smoothing of $A_1\cup {An}$. We require also that $\A \cap A_1'=\emptyset$. Set  $\A'=\{A_1',A_2,\dots,A_g\}$. Then $(\A',\B)$ is a new convex compressing disc system  and we say that $(\A,\B)$ and $(\A',\B)$ are related to each other by a \emph{triple point move}. See Figure \ref{fig:moveanst}.  

The analogous move for $\B$ is also called a triple point move.

\begin{proposition}\label{prop:abisot} [Triple point move] 
If $(\A',\B')$ is obtained from $(\A,\B)$ by a single triple point move, then   the refinements $\H(\A,\B)$ and $\H(\A',\B')$ have a common positive stabilisation.
\end{proposition}

\begin{proof}
Figure~\ref{fig:R3prep} shows the convex surfaces near the triple point itself.  Here, we assume that the entire move takes place within an $I$-invariant neighbourhood of  $\Sigma$. By assuming the isotopy is sufficiently small, we may ensure that the actual crossing remains away from the dividing sets on the compression discs, as shown in the figure.  
In the cases where one or both of the discs in the local model is a product disc, the associated refined Heegaard surfaces will be convexly isotopic and there is nothing to check.  We therefore turn our attention to the case when neither $A_i$ nor $B_j$ is a product disc, and curves $X_i$ and $X_j$ are shown in orange in Figure~\ref{fig:R3prep}.  In the refined splittings, a neighbourhood of each of these is added as a one-handle to the opposite handlebody.  The top picture in Figure~\ref{fig:move} shows the refined Heegaard surfaces.

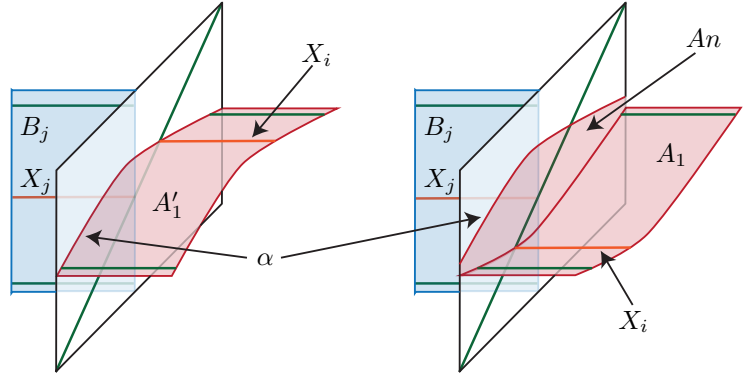
\begin{figure}[h]

\labellist
\small\hair 2pt
\pinlabel $An$ [l] at 272 137
\pinlabel $A_1$ [l] at 260 90
\pinlabel $\alpha$ [l] at 100 48
\pinlabel $A_1'$ [l] at 58 70
\pinlabel $B_j$ [l] at 5 100
\pinlabel $B_j$ [l] at 167 100
\pinlabel $X_i$ [l] at 118 130
\pinlabel $X_i$ [l] at 245 23
\pinlabel $X_j$ [l] at 5 79
\pinlabel $X_j$ [l] at 167 79

\endlabellist
\begin{center}
\includegraphics[scale=.95]{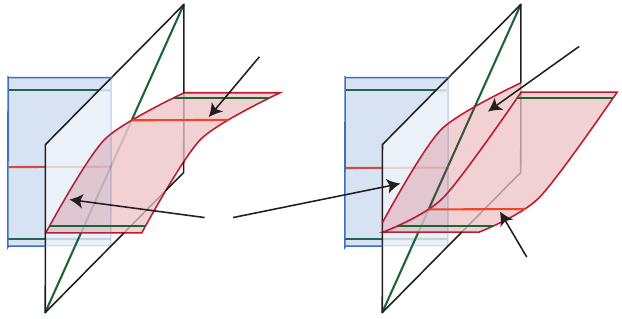}
\caption{ Convex surface model near a triple point.  The right-hand picture shows the original disc together with $An$ before smoothing, while the left hand picture shows the smoothed disc.} 
\label{fig:R3prep}
\end{center}
\end{figure}

In order to show that the two refinements have a common stabilisation, we identify a pair of discs, shown shaded in orange in the second row of Figure~\ref{fig:move}, which satisfy the hypotheses of Definition~\ref{def:posstab}. We attach $1$--handles to $V$ along the bold orange arcs in each picture and smooth the resulting Heegaard surfaces. It follows from Lemma~\ref{lem:stabHD} that the open books associated to the stabilised Heegaard diagrams are positive stabilisations of the originals.  On the other hand, one  may verify by inspection that the dividing sets on the two pictures in the bottom row are isotopic, as desired. The dividing curve is determined by 
the fact that $\tw_l(\xi,TD)=\frac{-1}{2}$: as seen in Figure~\ref{fig:move}, the dividing curve $\Gamma_{\partial\nu(l)}$  before smoothing rotates  half clockwise-turn less than $D\cap \partial \nu(l)$.

\begin{figure}[h]
\labellist
\small\hair 2pt
\pinlabel $V$ [l] at 10 270
\pinlabel $U$ [l] at 60 285
\pinlabel $V$ [l] at 190 285
\pinlabel $U$ [l] at 225 270

\pinlabel $l$ [l] at 100 165
\pinlabel $c$ [l] at 82 183
\pinlabel $l$ [l] at 270 217
\pinlabel $c$ [l] at 247 227

\endlabellist

\begin{center}
\includegraphics[scale=1]{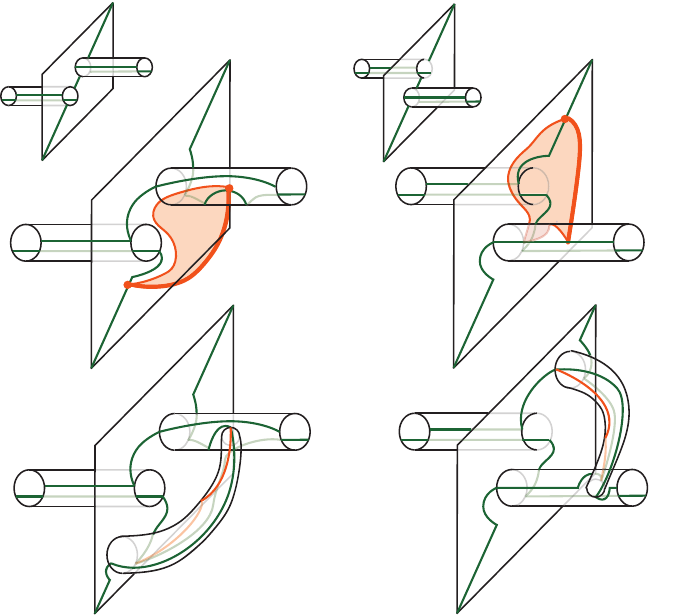}
\caption{ Inset: Heegaard surfaces associated to a triple point move, before smoothing the refining tunnels drilled along the curves $X_i$ on $\A$ and $X_j$ on $\B$ from Figure~\ref{fig:R3prep}.  Top: Add contact $1$--handles to $V$ along the bold $l$ curves.  The indicated $c$ curves show that these are positive stablisations. Bottom: After smoothing, the  stabilised  Heegaard surfaces are convexly isotopic. The orange curves in the final figure show the intersection of the orange discs with the added tunnels in order to make the relative twisting easier to see.} 
\label{fig:move}
\end{center}
\end{figure}
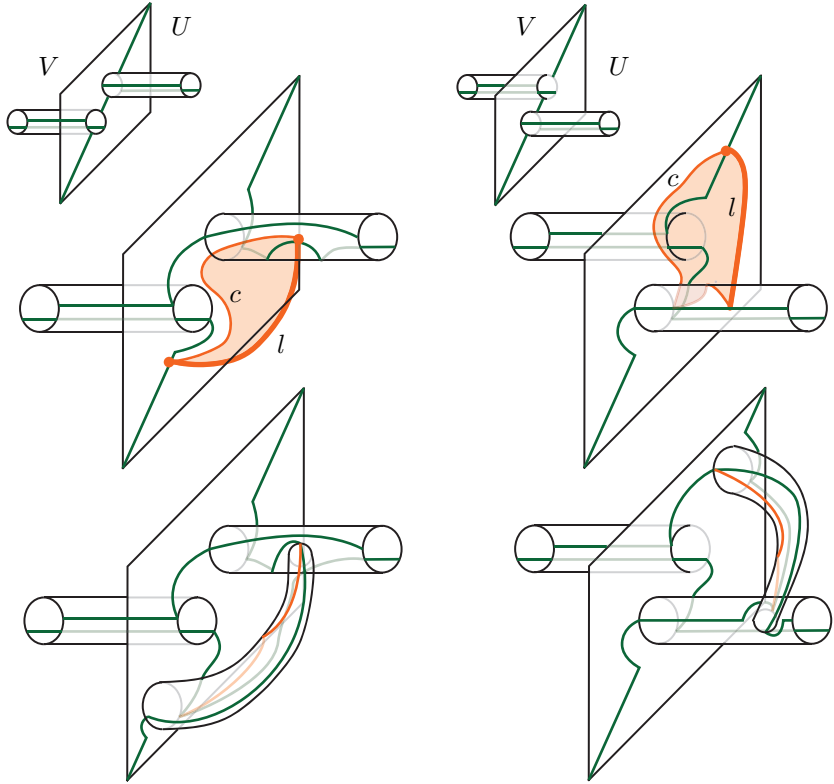

The proof in the case that ${\partial \B}$ crosses a point of ${\partial \A} \cap \Gamma_\Sigma$ is similar. 
\end{proof}

\subsubsection*{[F]: Finger move.} 

\begin{figure}[h]
\labellist
\small\hair 2pt
\pinlabel $An$ [l] at 38 78
\pinlabel $\alpha$ [l] at 118 17
\pinlabel $\partial A_1$ [l] at 138 42

\endlabellist

\begin{center}
\includegraphics[scale=1]{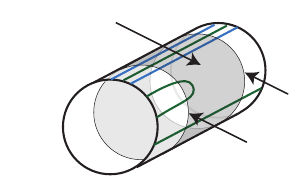}
\caption{ Topological model for the finger  [F] move on a convex compressing disc system.  } \label{fig:moveansf}
\end{center}
\end{figure}

Let $(\A,\B)$ be a convex compressing disc system and let $An\subset \Sigma$ be an annulus with Legendrian boundary $\partial A_1\cup \alpha$. Suppose that $\Gamma\cap An$ consists of parallel arcs from $A_1$ to $\alpha$, together with an additional boundary parallel arc  anchored at $\alpha$.
We further assume that $An$ is disjoint from $\partial \B\cap\Gamma$ and that $\partial \B\cap An$ consists only of parallel arcs between $\partial A_1$ and $\alpha$.
 Define $A_1'$ to be a convex surface properly embedded in  $U$ that is convex isotopic, relative to $\alpha$, to the smoothing of $A_1\cup An$.  We require also that $\A\cap A'_1 =\emptyset$. Set $\A'=\{A_1',A_2,\dots,A_g\}$.  Then  $(\A',\B)$  is a new convex compressing disc system and we say that $(\A,\B)$ and $(\A',\B)$ are related to each other by a \emph{finger move}. See Figure \ref{fig:moveansf}.

We use the term \emph{finger move} also to denote the reverse of this move, or the analogous moves of $\B$.

\begin{proposition}\label{prop:fingermove}[Finger move]\label{prop:finger}
If $(\A',\B')$ is obtained from $(\A,\B)$ by a single finger move, then   the refinements $\H(\A,\B)$ and $\H(\A',\B')$ have a common positive stabilisation.

\begin{figure}[h]
\labellist
\small\hair 2pt
\pinlabel $A_1$ [l] at 12 150
\pinlabel An [l] at 39 100
\pinlabel $A_1'$ [l] at 270 150

\endlabellist
\begin{center}
\includegraphics[scale=.9]{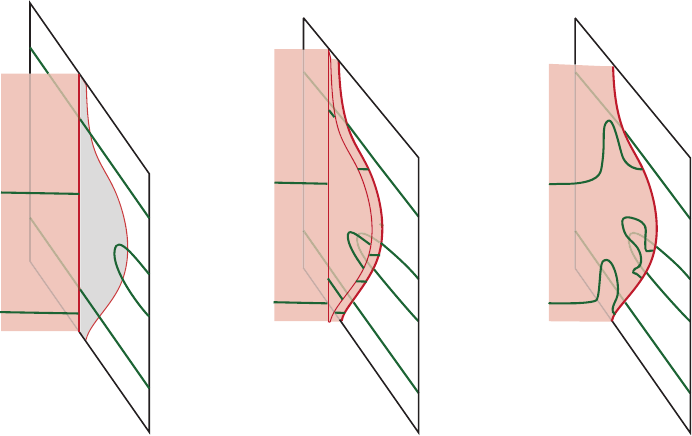}
\caption{ Left: The initial configuration for a finger move.  Centre: Use $An$ to build a a piecewise convex surface with Legendrian boundary. Right: A convex smoothing of the centre surface, showing the new dividing set on $A_1'$.}\label{fig:fingeran}\end{center}
\end{figure} 
\end{proposition}

\begin{proof} Consider Figures~\ref{fig:fingeran} and \ref{fig:fingeranstab}. The first picture in Figure~\ref{fig:fingeran} shows the original disc $A_1$ and the portion of the finger move annulus which contains the new intersection with $\Gamma_\Sigma$.  To construct $A_1'$,  smooth the piecewise convex surface shown in the central picture to get the new dividing set  shown on the right.  Since $\Gamma_{\A'}$ has an additional component, constructing the refinement requires tunneling along an additional arc $X'$, shown in orange in Figure~\ref{fig:fingeranstab}.  Observe that $X'$ cobounds a disc in $U$ with the orange arc on $\Sigma$ running parallel to the finger; these arcs are disjoint from $\Gamma$ except at their endpoints, so the twisting of $\xi$ relative to $TD$ is $-\frac{1}{2}$ in each case.    The hypotheses of Definition~\ref{def:posstab} are satisfied, so the open book associated to the new refinement is a positive stabilisation of the open book associated to the original refinement, as desired.
\begin{figure}[h]
\labellist
\small\hair 2pt

\pinlabel $A_1'$ [l] at 8 147
\pinlabel $X'$ [l] at 5 66

\endlabellist

\begin{center}
\includegraphics[scale=1]{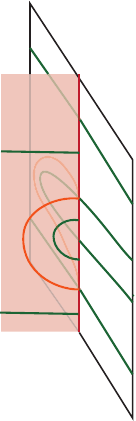}
\caption{  After performing a finger move that increases $|\Gamma_\Sigma\cap \partial \A|$ by two, refining the Heegaard splitting requires an additional positive stabilisation along $X'$.}\label{fig:fingeranstab}\end{center}
\end{figure}\end{proof}

\subsubsection*{[I]: Interior bypass}  
\begin{figure}[h]
\labellist
\small\hair 2pt
\pinlabel D [l] at 17 70
\pinlabel $A_1$ [l] at 38 90
\pinlabel $A_1'$ [l] at 88 8

\endlabellist

\begin{center}
\includegraphics[scale=1]{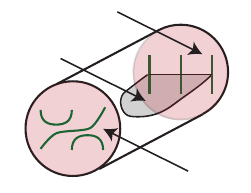}
\caption{ Local model  for the interior bypass  [I] move on a convex compressing disc system.  } \label{fig:moveansi}
\end{center}
\end{figure}

Let $(\A,\B)$ be a convex compressing disc system. Suppose that there is another convex compressing disc $A_1'$ for $U$ with Legendrian boundary that is obtained from $A_1$ by a bypass attachment along a bypass half-disc $D$ disjoint from $\A\setminus A_1$.
Then for $\A'=\{A_1',A_2,\dots,A_g\}$, we say that the convex compressing disc systems $(\A,\B)$ and $(\A',\B)$ are related to each other by an \emph{interior bypass move}. See Figure \ref{fig:moveansi}. 

\begin{proposition}\label{prop:intbyp}[Interior bypass]
If $(\A',\B')$ is obtained from $(\A,\B)$ by a single interior bypass move, then   the refinements $\H(\A,\B)$ and $\H(\A',\B')$ have   a  common positive stabilisation.
\end{proposition}

\begin{figure}[h]
 \labellist
\small\hair 2pt
\pinlabel $\Sigma_1'$ [l] at 5 30
\pinlabel $\Sigma_2'$ [l] at 228 30
\endlabellist

\begin{center}
\includegraphics[scale=.95]{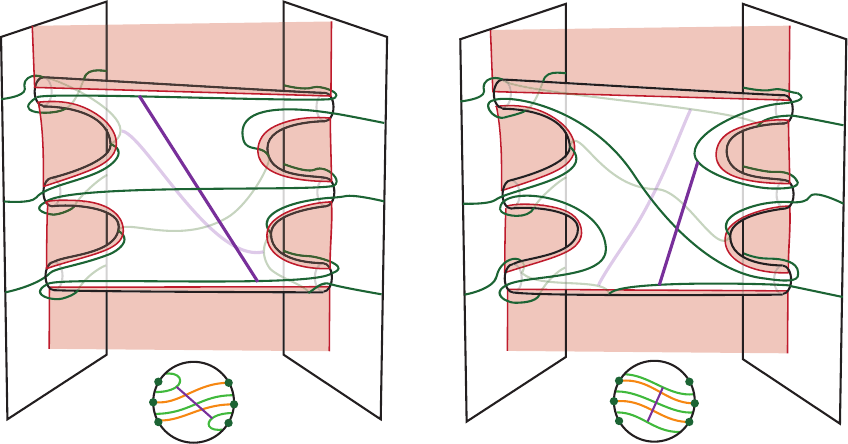}
\caption{The dividing sets on $\Sigma_1'$ and $\Sigma_2'$ become isotopic after performing a bypass along each of the indicated purple arcs.  Since the bypass attachment arcs lie in the complement of the indicated disc systems for the stabilised refinements, it follows from Proposition~\ref{prop:bypass} that the two convex Heegaard splittings admit a common positive stabilisation. } 
\label{fig:intbypassfin}
\end{center}
\end{figure}

\begin{proof}
Isotoping  the interior of a convex compressing disc $A_1$ across a bypass half-disc changes the dividing set $\Gamma_{A_1}$ as shown in the inset in Figure~\ref{fig:intbypass}.  Since we have previously shown that the refinement of $(\Sigma, \A, \B)$ is independent of the choice of arcs $\X$, we may choose $\X$ to agree with the orange curves shown in the inset picture.  Note also that these may be assumed to be local pictures and we place no restrictions on the other curves of $\Gamma_\A$ or $\X$; if the bypass arc $c$ intersects the same component of $\Gamma_\A$ twice, then the bypass is necessarily trivial and may be disregarded.

Figure~\ref{fig:intbypass} shows local models for the refined Heegaard surfaces $\Sigma_1$ and $\Sigma_2$ after tunneling along the chosen $\X$ before (left) and after (right) the isotopy across the bypass.  Each picture also shows three bold arcs labeled $l_i$ whose endpoints lie on the dividing set of the new tunnels.  Isotoping each $l_i$ arc to lie on the stabilised Heegaard surface traces out a disc whose boundary $l_i \cup c_i$ intersects $\Gamma_{\Sigma_j}$ only at the shared endpoints of the $l_i$ and $c_i$ arcs.  Since these discs lies in a neighbourhood of the original  $A_1$, their dividing sets are isotopic to the restriction of the original $\Gamma_{A_1}$, as shown. It follows that each of these discs satisfies the hypotheses of Definition~\ref{def:posstab},  so  Lemma~\ref{lem:stabHD} implies that stabilising the Heegaard splitting by tunneling along each $l_i$  arc stabilises the associated open book decomposition.  Performing these stabilisations yields the pair of isotopic Heegaard surfaces $\Sigma_1'$ and $\Sigma_2'$ shown at the bottom of Figure~\ref{fig:intbypass}, but the dividing sets $\Gamma_{\Sigma_1'}$ and $\Gamma_{\Sigma_2'} $ are not yet isotopic.

 \begin{figure}[h]
 
 \labellist
\small\hair 2pt
\pinlabel $l_1$ [l] at 62 384
\pinlabel $c_1$ [l] at 7 370

\pinlabel $l_2$ [l] at 69 334
\pinlabel $c_2$ [l] at 7 320

\pinlabel $l_3$ [l] at 124 296
\pinlabel $c_3$ [l] at 166 305

\pinlabel $l_1$ [l] at 250 313
\pinlabel $c_1$ [l] at 294 315

\pinlabel $l_2$ [l] at 360 389
\pinlabel $c_2$ [l] at 410 350

\pinlabel $l_3$ [l] at 342 344
\pinlabel $c_3$ [l] at 408 308

\pinlabel $\Sigma_1$ [l] at 10 250
\pinlabel $\Sigma_1'$ [l] at 10 30

\pinlabel $\Sigma_2$ [l] at 245 250
\pinlabel $\Sigma_2'$ [l] at 245 30

\endlabellist

\begin{center}
\includegraphics[scale=.9]{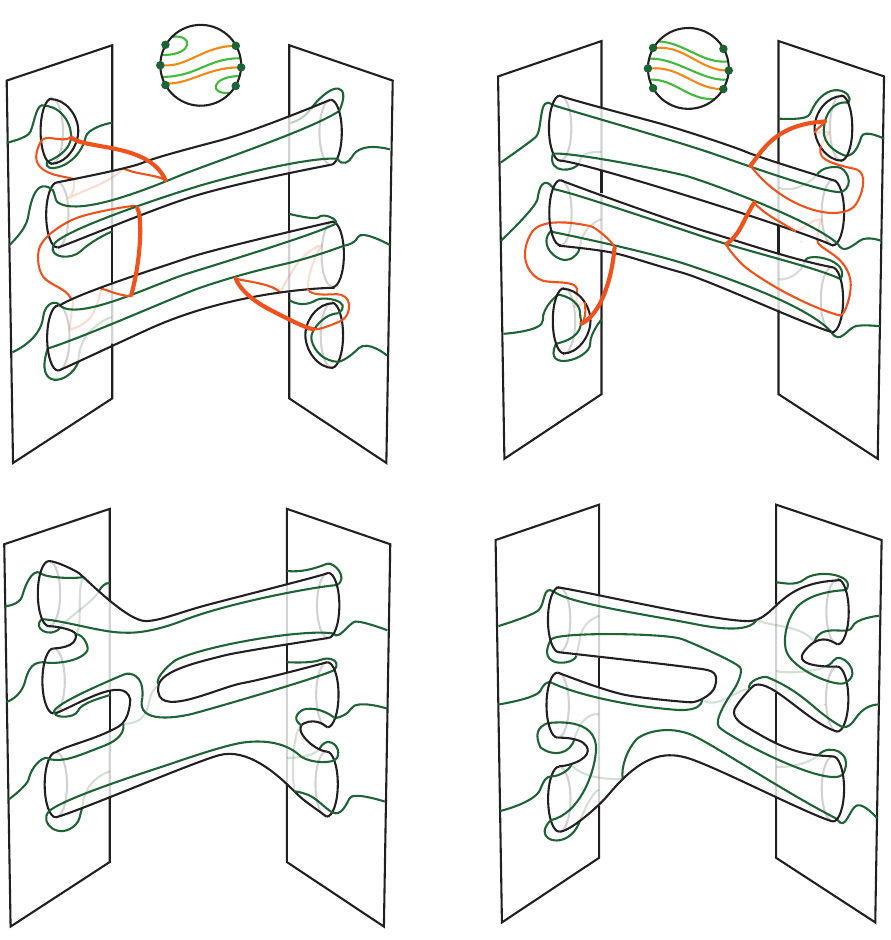}
\caption{When two compressing discs are related by a single bypass,  they give rise to distinct refinements $\Sigma_1$ and $\Sigma_2$.  Positively stabilising the Heegaard splittings by tunneling along the indicated bold $l_i$ arcs yields isotopic Heegaard surfaces $\Sigma_1 $ and $\Sigma_2$. } 
\label{fig:intbypass}
\end{center}
\end{figure}
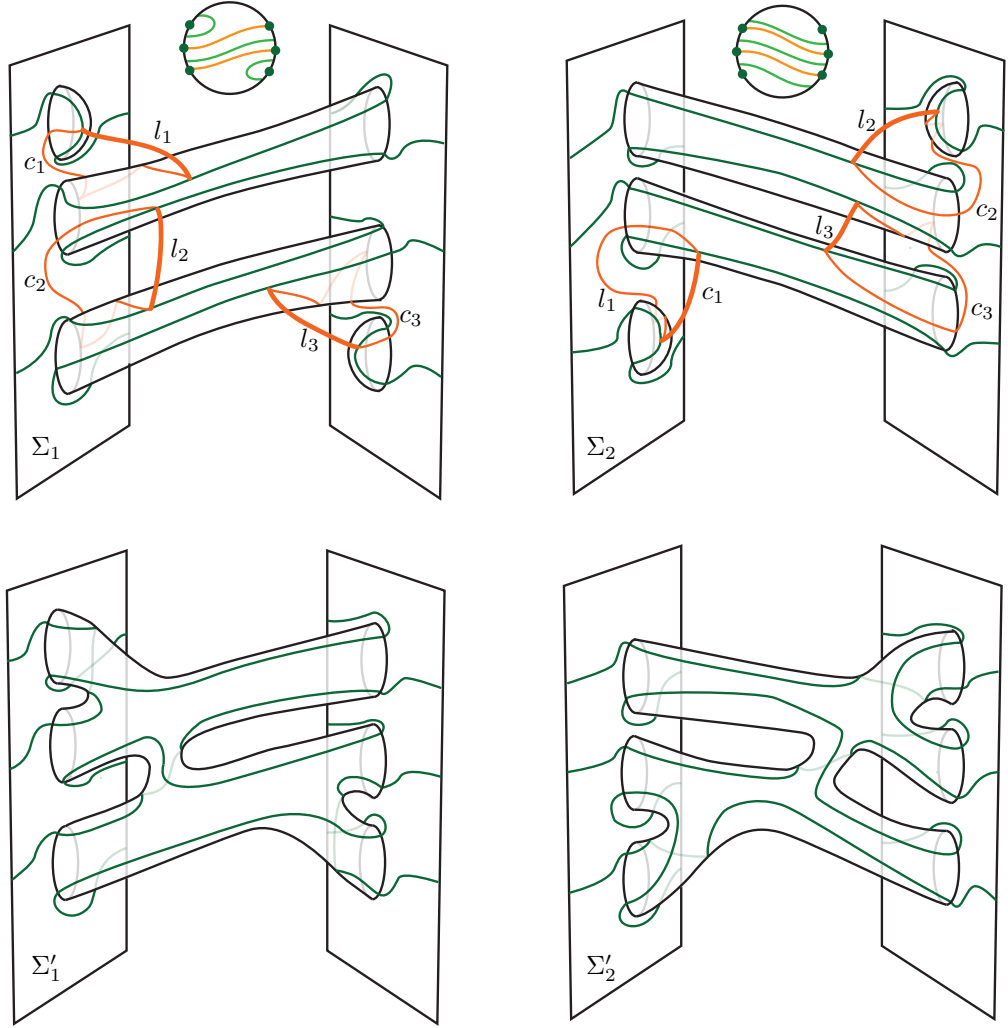

To complete the argument, we turn to Figure~\ref{fig:intbypassfin} where the isotopic Heegaard surface are shown with their dividing sets and the bypass attachment arcs that relate them.  Note that bypasses along these arcs necessarily exist, by the hypotheses of the interior bypass move.  Attaching a pair of bypasses  in the left-hand figure  renders the dividing set isotopic to that of the right-hand figure.  Note that with respect to the orientation of $\Sigma$ as $\partial U$, one of these bypasses is attached from the front and the other, from the back. Since the bypass arcs in each picture are disjoint from the indicated system of convex compressing discs for the stabilised refinements, Proposition~\ref{prop:bypass} implies that the two splittings admit a common positive stabilisation, as desired. 
\end{proof}

\subsubsection*{[H]: Handle slide.}
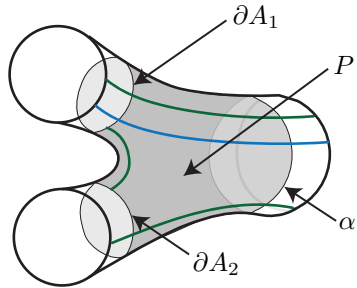
\begin{figure}[h]
\labellist
\small\hair 2pt
\pinlabel $P$ [l] at 138 85
\pinlabel $\partial A_1$ [l] at 100 103
\pinlabel $\partial A_2$ [l] at 85 10
\pinlabel $\alpha$ [l] at 140 25

\endlabellist

\begin{center}
\includegraphics[scale=1]{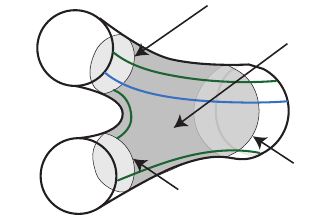}
\caption{ Local model for the handleslide [H] move on a convex compressing disc system.  } \label{fig:moveans2h}
\end{center}
\end{figure}

Let $(\A,\B)$ be a convex compressing disc system and consider a pair of pants  $P\subset \Sigma$ with Legendrian boundary $\partial A_1\cup \partial A_2 \cup \alpha$. Suppose that $\Gamma\cap P$ consists of a single arc $\gamma$ connecting $\partial A_1$ to $\partial A_2$ and an arbitrary number of arcs connecting $\partial A_1\cup \partial A_2$ to $\alpha$.
Suppose also that $P$ is disjoint from $\partial \B\cap\Gamma$ and that $\partial \B\cap P$  consists only of {disjoint} arcs connecting  $\partial A_1\cup \partial A_2$ to $\alpha$; each of these arcs is required to be disjoint from $\gamma$.  See Figure~\ref{fig:moveans2h}.

 Let $A_1'$ be a convex surface properly embedded in $U$ that is convex isotopic, relative to $\alpha$, to the smoothing of $P\cup A_1\cup A_2$.  We require also that $A_1'$ is  disjoint from $\A$. Then for $\A'=\{A_1',A_2,\dots,A_g\}$, we get a new convex compressing disc system $(\A',\B)$. We say that $(\A,\B)$ and $(\A',\B)$ are related to each other by a \emph{handle slide along $\gamma$}. See Figure \ref{fig:moveans2h}.  
We can similarly perform handle slides across $\B$ discs.  

\begin{proposition}\label{prop:handle slide}[Handle slide]
If $(\A',\B')$ is obtained from $(\A,\B)$ by a single handle slide, then   the refinement of $\H$ via $\A$ and $\B$ and  the refinement of $\H$ via $\A'$ and $\B'$ have a common positive stabilisation.
\end{proposition}

\begin{proof}

As a first step, we observe that any collection of convex discs with Legendrian boundary on $\Sigma$ defines a refinement as long as there is a subset of the discs which constitute a disc system for each handlebody.  The extension is immediate:  choose sufficiently many $\X$ arcs to separate components of $\Gamma_{\A\cup \B}$ in every disc and add tunnels along all of these.

We will show that the refinement of $(\{A_1, A_1', A_2, \dots, A_g\}, \B)$ is a positive stabilisation of each of $(\A, \B)$ and $(\A', \B)$.  In fact, it suffices to show that $(\{A_1, A_1', A_2, \dots, A_g\}, \B)$ is a positive stabilisation of $(\A, \B)$, as the relationship between the curves $\{A_1, A_2, A_1'\}$ is symmetric up to finger moves which have already been shown to preserve the positive stabilisation class.

Given the pair of pants defining the handle slide move, we can explicitly construct $A_1'$ by smoothing a piecewise convex surface built from parallel copies of $A_1$ and $A_2$ and neighbourhood of $\gamma$ on $\Sigma$.
Observe in Figure~\ref{fig:handleslide} that the dividing set and choice of $\X$ on $A_1$ and $A_2$ dictate the dividing set and a canonical choice of $\X$ on $A_1'$.  The $\X$ arcs on $A_1'$ come in two forms: arcs that are parallel to $\X$ arcs on $A_1$ or $A_2$ and arcs which cross the band defined by $\gamma$.  

 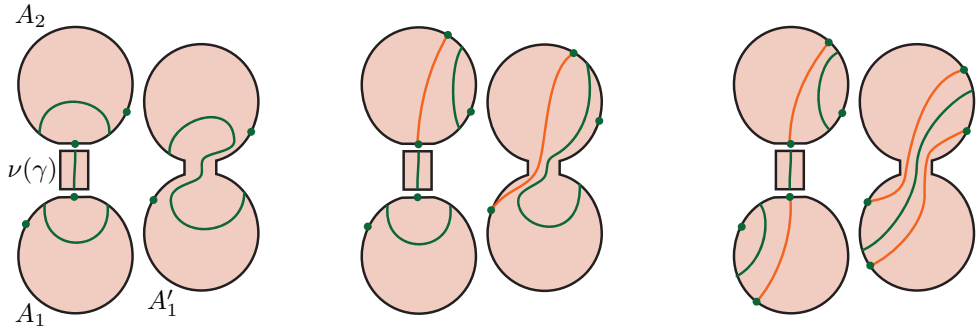
\begin{figure}[h]
 \labellist
\small\hair 2pt
\pinlabel $A_1$ [l] at -3 0
\pinlabel $A_1'$ [l] at 50 5
\pinlabel $A_2$ [l] at -3 120
\pinlabel $\nu(\gamma)$ [l] at -6 58

\endlabellist

\begin{center}
\includegraphics[scale=.95]{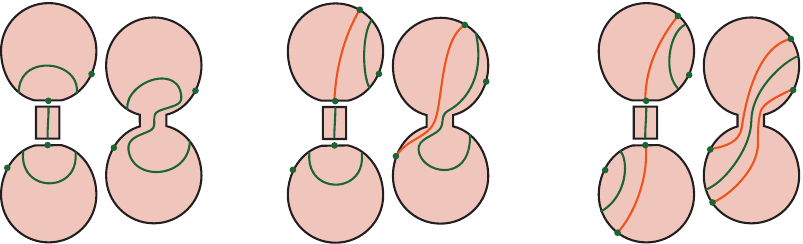}
\caption{Local models for $A_1'$.  The orange curves in the glued-up discs are arcs of $\X$ in $\A_i'$ which cross the band $\nu(\gamma)$.} 
\label{fig:handleslide}
\end{center}
\end{figure}

We consider the parallel arcs first.  Let $X_1\subset A_1$ be an arc which is tunneled along in the refinement via $\A$ and  $\B$.  We see in Figure~\ref{fig:parallelhandle} that tunneling along a copy of the same arc in a disc locally parallel to $A_1$ is a positive stabilisation.

\begin{figure}[h]
 \labellist
\small\hair 2pt
\pinlabel $A_1$ [l] at -3 151
\pinlabel $X_1$ [l] at 30 134

\pinlabel $l$ [l] at 95 55
\pinlabel $c$ [l] at 156 76

\endlabellist

\begin{center}
\includegraphics[scale=.95]{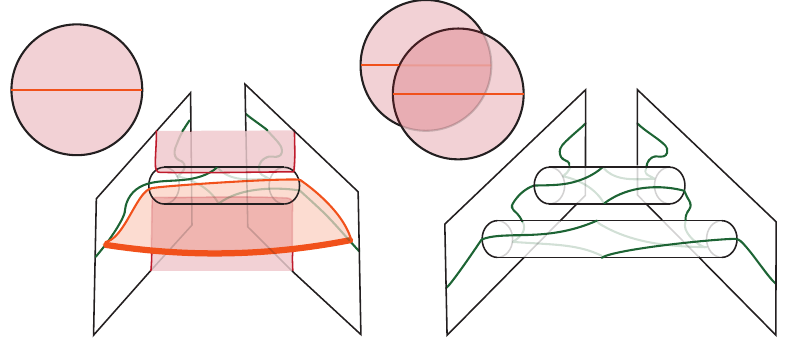}
\caption{ Left inset: A copy of $A_1$ with an $X_1$ arc shown in orange.  Right inset: Parallel copies of $A_1$ with parallel copies of $X_1$. Left: After adding a tunnel along $X_1$, the  shaded disc with boundary $l\cup c$ satisfies the hypotheses of Definition~\ref{def:posstab}.  The bold arc $l$ is the second copy of $X_1$.  Right: Adding a tunnel along $l$ is a positive stabilisation of the original refinement. } 
\label{fig:parallelhandle}
\end{center}
\end{figure}
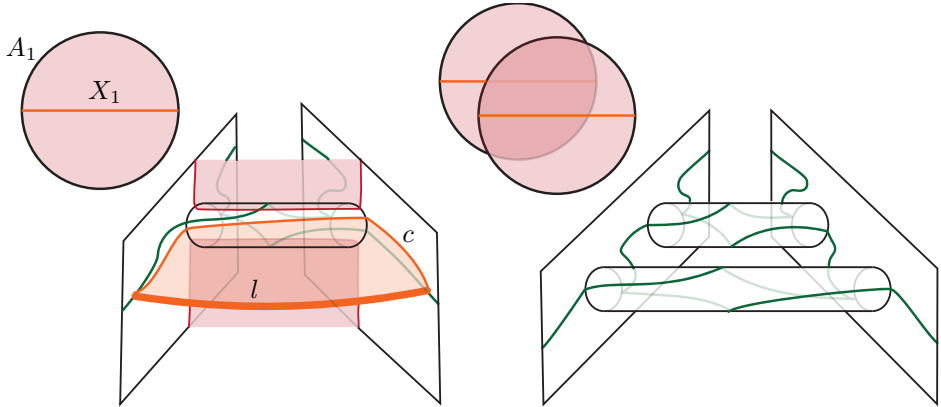

Finally, we consider the $\X$ arcs that start cross the band defined by $\gamma$, as shown in Figure~\ref{fig:handleslide}.  Figure~\ref{fig:pantstunnel} shows that tunneling along each of these arcs is a positive stabilisation of the refinement via the original $ \A$ and $\B$, as desired.
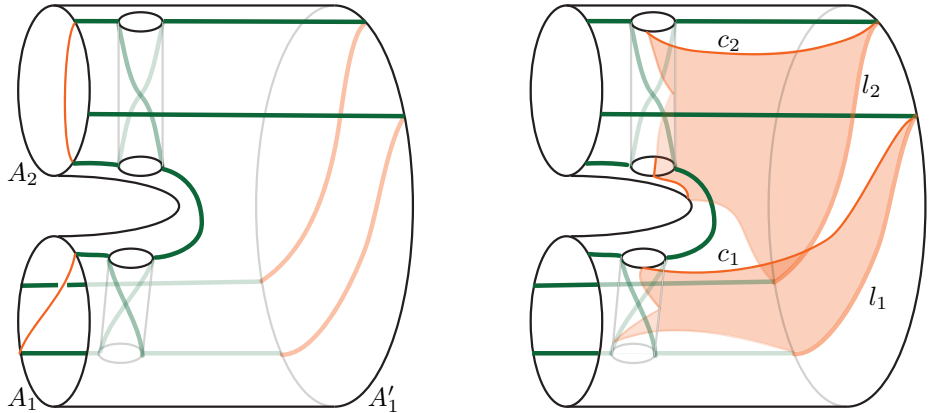
\begin{figure}[h]
 \labellist
\small\hair 2pt
\pinlabel $A_1$ [l] at -5 5
\pinlabel $A_2$ [l] at -5 105
\pinlabel $A_1'$ [l] at 155 5
\pinlabel $c_1$ [l] at 310 68
\pinlabel $l_1$ [l] at 377 50
\pinlabel $c_2$ [l] at 310 163
\pinlabel $l_2$ [l] at 374 144
\endlabellist

\begin{center}
\includegraphics[scale=.85]{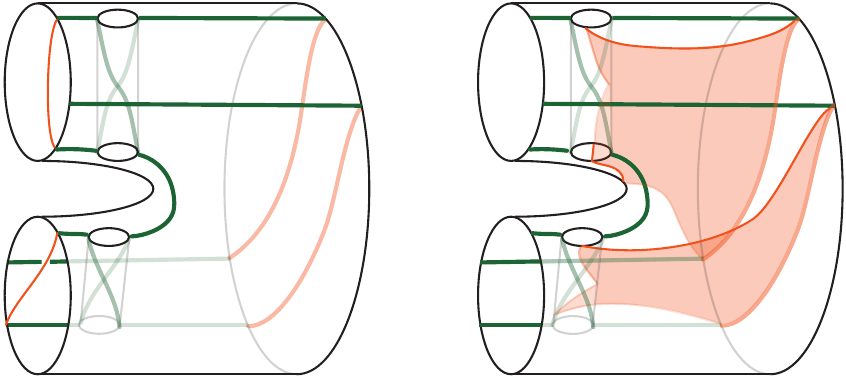}
\caption{  Left: Orange $\X$ arcs on $A_1$ and $A_2$ guide the tunnels shown; following Figure~\ref{fig:handleslide}, the bold arcs indicate where tunnels should be excavated on  $A_1'$.  Right: The shaded discs bounded by $l_i \cup c_i$ show that tunneling along the bold $l_i$ arcs  positively stabilises the original splitting.} 
\label{fig:pantstunnel}
\end{center}
\end{figure}
\end{proof}

\subsubsection{Sufficiency of disc moves}\label{sec:suff}
We conclude this section, and the proof of Theorem~\ref{thm:indep}, by showing that the moves introduced above act transitively on the equivalence classes of  convex compressing disc systems for a tight Heegaard splitting. 
\begin{proposition}\label{prop:elmoves}
Up to convex isotopy, any two convex compressing discs systems $(\A,\B)$ and $(\A',\B')$ for a given tight Heegaard splitting $\H$ of $(M,\xi)$ are related by a sequence of  elementary  disc moves \textit{[F]},  \textit{[T]}, \textit{[I]} and \textit{[H]}.  
\end{proposition}
\begin{proof}
If $(\A,\B)$ and $(\A',\B')$ are smoothly isotopic, then we first multiply the isotopy by a cut-off function that vanishes away from a neighbourhood $N'$ of $N$ that is an extension  of $N$ by $I$-invariant neighbourhoods of $\Sigma$ and $\Sigma'$.
 This isotopy can be followed by finger moves between $\partial \A$ and $\partial \B$ which have no effect on the disc system and the elementary moves \textit{[F]} and \ \textit{[T]}.   To be able to do this, at each step one needs to simultaneously Legendrian realise $\partial \A\cup \partial A_1'$, but as $\partial A_i\cap\Gamma\neq\emptyset$, this is always possible by a $C^\infty$-isotopy of $\Sigma$. 

Now fixing the boundary of the compressing discs,  Theorem \ref{thm:isdisc} decomposes the rest of the isotopy as a composition of bypasses from the front or from the back. This is move  \textit{[I]}. 

Smoothly, any two systems of compressing discs can be made isotopic by handle slides; by first performing some isotopies and finger moves in $N'$, one can make sure that the handle slide  is performed along an arc of $\Gamma$, as required by the local model for move \textit{[H]}.  
This finishes the proof.  
\end{proof}
We have now established everything needed to conclude the invariance of the positive stabilisation class of the open book associated to a tight splitting:

\begin{proof}[Proof of Theorem \ref{thm:indep}.]
By Proposition~\ref{prop:elmoves} the two convex compressing disc systems are related to each other by elementary disc moves, and  Propositions~\ref{prop:abisot}, \ref{prop:fingermove}, \ref{prop:intbyp}, and \ref{prop:handle slide} show that each of these elementary disc moves  preserves the positive stabilisation class of the open book associated to the refinement. 
\end{proof}

\section{Proof of the Giroux Correspondence}\label{sec:GC} 

We now have all the ingredients to prove the Giroux Correspondence for tight contact $3$--manifolds.

 \begin{theorem}[c.f.\ Theorem \ref{thm:GC}]\label{thm:GC2} Let $(M, \xi)$ be a tight contact manifold and suppose that $(B, \pi)$ and $(B',\pi')$ are two open book decompositions of $M$ supporting $\xi$.  Then $(B, \pi)$ and $(B',\pi')$ admit a common positive stabilisation. 
   \end{theorem}

  \begin{proof}
Construct the convex Heegaard splittings $\H=\H(B,\pi)$ and $\H'=\H(B',\pi')$ corresponding to $(B,\pi)$ and $(B',\pi')$, respectively. 

By the Reidemeister-Singer Theorem, $\H$ and $\H'$ will become isotopic after sufficiently  many Heegaard splitting stabilisations.  This topological statement places no restrictions on the stabilisations, so in the case where the Heegaard surfaces are convex in a contact manifold, we may choose each stabilisation to be positive in the sense of Definition~\ref{def:posstab}.
Stabilising $\H$ and $\H'$ accordingly results in a pair of convex Heegaard splitttings  of $(M,\xi)$ that are smoothly isotopic. By Lemma~\ref{lem:stabHD}, it then suffices to prove the theorem for a pair of open books that induce smoothly isotopic convex Heegaard splittings. 

Suppose now that the convex Heegaard splittings $\H=\H(B,\pi)=(\Sigma, U, V)$ and $\H'=\H(B',\pi')=(\Sigma', U',V')$ are smoothly isotopic. As noted in Theorem~\ref{thm:isdisc}, isotopy discretisation implies that $\Sigma$ and $\Sigma'$ are related by a sequence of bypasses, so we may enumerate the intermediate convex splitting surfaces $\Sigma=\Sigma_0,\Sigma_1,\dots,\Sigma_k=\Sigma'$.

We claim that for any consecutive pair of Heegaard surfaces $\Sigma_i, \Sigma_{i+1}$ defining tight Heegaard splittings, there exists a convex compressing disc system disjoint from the bypass half-disc.   By Theorem~\ref{thm:indep},  the positive stabilisation class of the splitting is independent of the choice of the convex compressing disc system, so it  follows from Corollary~\ref{cor:refstab} that the refinements of the splittings admit a common positive stabilisation.  Since this holds for each  pair in the sequence, we conclude that $\H$ and $\H'$ admit a common positive stabilisation.  

Applying Lemma~\ref{lem:stabHD} yet again, it follows that $(B,\pi)$ and $(B', \pi')$ admit a common positive open book stabilisation, as desired.

To complete the argument, we need only prove the claim that a single bypass attachment may always be performed in the complement of a convex compressing disc system. 

Choose a convex compressing disc system $(\A, \B)$ for the tight Heegaard splitting $(\Sigma, U, V)$.    Perform finger moves on $\partial \A$ and $\partial \B$ along the Legendrian attaching arc $c$ of the bypass half-disc $D$ until the boundaries of the compressing discs are all disjoint from $c$.  Now consider  intersections between $\B$ and $D$.  If  $\B\cap D$ has any simple closed curve components, we may isotope the interior of $\B$ across $l\subset \partial D$  until $\B\cap D$ consists only of arcs properly embedded in $D$.  Starting from outermost arcs, perform internal bypasses  to push $\B$ off $D$.  The result is a new system of convex compressing discs that is disjoint from $D$, as desired.  
\end{proof}

\section{Overtwisted manifolds} 

In this final section, we consider the case of overtwisted manifolds.  

With the exception of Section~\ref{sec:GC}, we have focused on tight Heegaard splittings rather than tight manifolds, so most of the technical results apply equally well to overtwisted contact manifolds divided into two tight handlebodies.  For example,  the last four paragraphs of the proof of Theorem~\ref{thm:GC}, together with Corollary~\ref{cor:refstab}, establish the following lemma for an arbitrary contact $3$-manifold:

\begin{lemma} Suppose that $(M,\xi)$ is a contact manifold with tight Heegaard splittings $\H=(\Sigma, U, V)$ and $\H'=(\Sigma, U', V')$ such that the handlebody $U'$ is built from $U$ by attaching  a single bypass slice. Let $\H''$ be the Heegaard splitting formed by attaching to $U$ only the $1$-handle associated to this bypass.  Then the refinements of $\H$, $\H'$, and $\H''$ admit common positive stabilisations.
\end{lemma}

In fact, we may generalise the Heegaard splittings we consider yet further by demanding more of the convex compressing disc systems.

\begin{definition} A \emph{tightening system} is a convex compressing disc system for $(\Sigma, U, V)$ with the property that each of $U\setminus \A$ and $V\setminus \B$ is a tight ball.
\end{definition} 

When $(\Sigma, U, V)$ is a tight Heegaard splitting, then every convex compressing disc system is a tightening system.  However, the requirement that the handlebodies $U$ and $V$ be tight is stronger than needed in order to define the refinement of the splitting.  

\begin{lemma} Let $(\A, \B)$ be a tightening system for the splitting $\H=(\Sigma, U, V)$.  Then the refinement of  $\H$ is a convex Heegaard splitting.
\end{lemma}

The proof of Proposition~\ref{prop:virtcont} began by cutting $U$ along $\A$ and $V$ along $\B$ to get a pair of tight balls; since this is the defining property of a tightening system, the argument holds as written in the case when the original splitting was tight. This suggests a path towards proving the full Giroux Correspondence: if we can arrange a tightening system at each step, then our sequence of Heegaard splittings will again produce a sequence of open books in a fixed positive stabilisation class. The final step in this program is ensuring that tightening systems may be found in the complement of bypass discs, and here we are only partially successful.  

We begin with the good news.  A potential obstruction to passing $\Sigma$ across a bypass half-disc in $V$ occurs if a given tightening system $\A \subset U$ intersects the bypass attachment arc on $\Sigma$.  However, a careful analysis of internal bypasses allows us to replace any such $\A$ with a different tightening system  for $U$ which remains a tightening system after the bypass attachment.  

Unfortunately, this is only half the battle, as there are Heegaard splittings without tightening systems.  The crux of the difficulty is the well known fact that overtwisted discs come in families; although it is trivial to choose a disc system that intersects a fixed overtwisted disc in an essential way, it is not always possible to kill all the overtwisted discs by intersection.
 
 \begin{example}[Heegaard splitting that does not admit a tightening system] Consider any genus-one convex Heegaard surface where the dividing set on the solid torus $U$ is a pair of parallel meridians.  It follows that  any compressing disc  has only inessential intersections with $\Gamma$.  If $\Gamma\cap \partial \A=\emptyset$, then any Legendrian realisation of $\partial \A$ is a $\tb=0$ unknot, and hence bounds an overtwisted disc.  If $|\Gamma\cap \partial \A|>2$, the proof of Proposition~\ref{prop:finger}  shows that reducing the number of intersections via a finger move preserves the isotopy class of the dividing set on the smoothed ball $U\setminus \A$.  We may thus assume $|\Gamma\cap \partial \A|=2$.   In this case,  cutting along $\A$ and smoothing yields a ball with a disconnected dividing set, so the ball remains overtwisted. This shows that there is no tightening system for such a splitting.
 \end{example}

The final example shows that tightening systems do sometimes exist.  

 \begin{example}[Overtwisted handlebody with a tightening system]  Consider a solid torus $U$ with a convex meridional disc $A$ with Legendrian boundary.  The disc $A$ is a tightening system for the solid torus if, after cutting along $A$ and smoothing, the dividing set on the ball is connected.  The left-hand side of Figure~\ref{fig:tightening} shows a connected dividing set on the boundary of a ball, where the surface is decomposed as an annulus together with two red discs indicating where the cut along $A$ occurred.  The figure also shows two copies of an arc on $A$ indicating where a bypass half-disc meets $A$ in the solid torus. Note first that the bypass from the back along the dashed purple arc is trivial, so it necessarily exists in $U\setminus A$.
  Passing $A$ across the bypass attaches a half-disc to the orange curve from the front and to the purple arc from the back.  Denote the new meridional disc by $A'$.   The right-hand figure shows the dividing set on the ball formed by cutting the solid torus along $A'$; since the dividing set is disconnected, the ball is overtwisted, showing  that the original contact structure on $U$ was overtwisted.
 
  \begin{figure}[h]
\begin{center}
\includegraphics[scale=1.3]{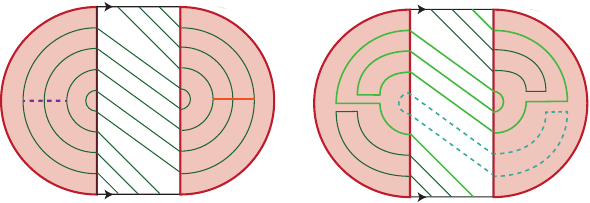}
\caption{ Left: Identifying the red discs and the two black intervals reconstructs a torus. The dashed purple and solid orange arcs indicate the location of a bypass half-disc in the solid torus.  Right: Passing the meridional disc across the bypass yields a disconnected dividing set.  } 
\label{fig:tightening}
\end{center}
\end{figure}
\end{example}

\bibliographystyle{alpha}

\bibliography{fob}

\end{document}